\documentstyle[12pt]{amsart}

\title[Positive Scalar Curvature]{Positive Scalar Curvature and Minimal Hypersurface Singularities}
\author{Richard Schoen}
\address{Department of Mathematics \\
                 University of California, Irvine \\
                 Irvine, CA 92697}

\thanks{The first author was partially supported by NSF grant
DMS-1404966.}

\author{Shing-Tung Yau}
\address{Department of Mathematics \\
Harvard University \\
Cambridge, MA 02138}

\thanks{The second author was partially supported by NSF grants DMS-1308244
and PHY-1306313}

\date{\today}

\newtheorem{thm}{Theorem}[section]
\newtheorem{lem}[thm]{Lemma}

\newtheorem{prop}[thm]{Proposition}

\theoremstyle{definition}
\newtheorem{rem}{Remark}[section]

\newtheorem{defn}{Definition}[section]

\numberwithin{equation}{section}

\renewcommand{\a}{\alpha}

\renewcommand{\d}{\delta}

\newcommand{\e}{\varepsilon}

\newcommand{\p}{\partial}

\newcommand{\s}{\sigma}
\newcommand{\Sig}{\Sigma}

\renewcommand{\O}{\Omega}
\renewcommand{\o}{\omega}

\newcommand{\z}{\zeta}

\def\Pb{\ifmmode{\Bbb P}\else{$\Bbb P$}\fi}
\def\Z{\ifmmode{\Bbb Z}\else{$\Bbb Z$}\fi}
\def\Q{\ifmmode{\Bbb Q}\else{$\Bbb Q$}\fi}
\def\C{\ifmmode{\Bbb C}\else{$\Bbb C$}\fi}
\def\R{\ifmmode{\Bbb R}\else{$\Bbb R$}\fi}
\def\S{\ifmmode{S^2}\else{$S^2$}\fi}

\def\S{\cal S}

\begin{document}

\maketitle

\begin{abstract} In this paper we develop methods to extend the minimal hypersurface
approach to positive scalar curvature problems to all dimensions. This includes a proof
of the positive mass theorem in all dimensions without a spin assumption. It also includes
statements about the structure of compact manifolds of positive scalar curvature extending
the work of \cite{sy1} to all dimensions. The technical
work in this paper is to construct minimal slicings and associated weight functions in the 
presence of small singular sets and to show that the singular sets do not become too 
large in the lower dimensional slices. It is shown that the singular set in any slice is a 
closed set with Hausdorff codimension at least three. In particular for arguments which 
involve slicing down to dimension $1$ or $2$ the method is successful. The arguments
can be viewed as an extension of the minimal hypersurface regularity theory to this setting of
minimal slicings.
\end{abstract}

\setcounter{secnumdepth}{1}

\setcounter{section}{0}

\section{\bf Introduction}

The study of manifolds of positive scalar curvature has a long history in
both differential geometry and general relativity. The theorems involved 
include the positive mass theorem, the topological classification of manifolds 
of positive scalar curvature, and the local geometric study of metrics of
positive scalar curvature. There are two methods which
have been successful in this study in general situations, the Dirac operator
method and the minimal hypersurface method. Both of these methods have
restrictions on their applicability, the Dirac operator methods require the 
topological assumption that the manifold be spin, and the minimal hypersurface
method has been restricted to the case of manifolds with dimension at most $8$
because of the possibility of singularities which might occur in the hypersurfaces.
The purpose of this paper is to extend the minimal hypersurface method to
all dimensions. 

The Dirac operator method was pioneered by A. Lichnerowicz \cite{lich} and
M. Atiyah, I. Singer \cite{as} in the early 1960s. It was extended by N. Hitchin \cite{h}
and then systematically developed by M. Gromov and H. B. Lawson in \cite{gl1}, \cite{gl2}, and
\cite{gl3}. Surgery methods for manifolds of positive scalar curvature were developed 
in \cite{sy1} and \cite{gl2}. For simply connected manifolds $M^n$ with $n\geq 5$ Gromov 
and Lawson conjectured necessary and conditions for $M$ to have a metric of
positive scalar curvature (related to the index of the Dirac operator in the spin case).
The conjecture was solved in the affirmative by S. Stolz \cite{st}. The Dirac operator method
was used by E. Witten \cite{w} to prove the positive mass theorem for spin manifolds (see
also \cite{pt}).

The minimal hypersurface method originated in \cite{sy4} for the three dimensional case and
was extended to higher dimensions in \cite{sy1}. The extension to the positive mass theorem was initiated in \cite{sy2} and in higher dimensions in \cite{sy5} and \cite{sc}. In this paper we extend
the minimal hypersurface argument to all dimensions at least as regards the applications
to the positive mass theorem and results which can be proven by slicing down to dimension
two. 

The basic objects of study in this paper are called {\it minimal $k$-slicings} and we now
describe them. We start with a compact oriented Riemannian manifold $M$ which will be our
top dimensional slice $\Sig_n$. We choose an oriented volume minimizing hypersurface
$\Sig_{n-1}$. Since $\Sig_{n-1}$ is stable, the second variation form $S_{n-1}(\varphi,\varphi)$
has first eigenvalue which is non-negative. We choose a positive first eigenfunction $u_{n-1}$
and we use it as a weight $\rho_{n-1}$ for the volume functional on $n-2$ cycles which are
contained in $\Sig_{n-1}$. We assume we have a $\Sig_{n-2}\subset\Sig_{n-1}$ which
minimizes the weighted volume $V_{\rho_{n-1}}(\cdot)$. The second variation $S_{n-2}(\varphi,\varphi)$ for the weighted volume on $\Sig_{n-2}$ then has non-negative first eigenvalue and
we let $u_{n-2}$ be a positive first eigenfunction. We then define $\rho_{n-2}=u_{n-2}\rho_{n-1}$
and we continue this process. That is if we have $\Sig_{j+1}\subset\Sig_{j+2}\subset\ldots\subset
\Sig_n$ which have been constructed, we choose $\Sig_j$ to be a minimizer of the weighted
volume $V_{\rho_{j+1}}(\cdot)$. Such a nested family $\Sig_k\subset\Sig_{k+1}\subset\ldots\subset
\Sig_n$ is called a {\it minimal $k$-slicing}. 

The basic geometric theorem about minimal $k$-slicings which is generalized in Section 2
is the statement that if $\Sig_n$ has positive scalar curvature then for any minimal $k$-slicing
we have that $\Sig_k$ is Yamabe positive and so admits a metric of positive scalar curvature.
In particular if $k=2$ then $\Sig_2$ must be diffeomorphic to $S^2$ and there can be no
minimal $1$-slicing.

If we start with $\Sig_n$ with $n\geq 8$, there might be a closed singular set ${\cal S}_{n-1}$ of Hausdorff dimension at most $n-8$ in $\Sig_{n-1}$. In this paper we develop methods to
carry out the construction of minimal $k$-slicings allowing for the possibility that the $\Sig_j$
may have nonempty singular sets ${\cal S}_j$. In order to do this it is necessary to extend
the existence and regularity theory for minimal hypersurfaces to this setting. To do this requires
maintaining some integral control of the geometry of the $\Sig_j$ in the ambient manifold $\Sig_n$,
and also of constructing the eigenfunctions $u_j$ which are bounded in appropriate weighted
Sobolev spaces. This control is gotten by carefully exploiting the terms which are left over
in the geometry of the second variation at each stage of the slicing. This is done by modifying
the second variation form $S_j$ to a larger form $Q_j$. The form $Q_j$ is more coercive
and can be diagonalized with respect to the weighted $L^2$ norm even in the presence of
small singular sets. We can then construct the next slice using the first eigenfunction for the
form $Q_j$ to modify the weight. This procedure only works if the singular sets ${\cal S}_j$
do not become too large. We prove that for a minimal $k$-slicing the Hausdorff dimension
of the singular set ${\cal S}_k$ is at most $k-3$. The regularity theorem is proven by
establishing appropriate compactness theorems for minimal $k$-slicings and showing that
at a singular point there is a homogeneous minimal $k$-slicing gotten by rescaling and
using appropriate monotonicity theorems (volume monotonicity and monotonicity of an
appropriate frequency function). A homogeneous minimal $k$-slicing is one in ${\mathbb R}^n$
for which all of the $\Sig_j$ are cones and all of the $u_j$ are homogeneous of some degree.
It is then possible to show that if we had a $\Sig_{k+1}$ with singular set of codimension
at least $3$, but $\Sig_k$ had a singular set of Hausdorff dimension larger then $k-3$ then there
would exist a nontrivial homogeneous $2$-slicing with $\Sig_2$ having an isolated singularity
at the origin. We show that no such homogeneous slicings exist to conclude that if ${\cal S}_{k+1}$
has codimension at least $3$ in $\Sig_{k+1}$, then ${\cal S}_k$ has codimension at least $3$
in $\Sig_k$. In particular if $k=2$ then $\Sig_2$ is regular. 

We now state the main theorems of the paper beginning with the positive mass theorem.
A manifold $M^n$ is called asymptotically flat if there is a compact set $K\subset M$
such that $M\setminus K$ is diffeomorphic to the exterior of a ball in ${\mathbb R}^n$
and there are coordinates near infinity $x^1,\ldots, x^n$ so that the metric components
$g_{ij}$ satisfy
\[ g_{ij}=\d_{ij}+O(|x|^{-p}),\ |x||\partial g_{ij}|+|x|^2|\partial^2g_{ij}|=O(|x|^{-p}) 
\]
for some $p>\frac{n-2}{2}$. We also require the scalar curvature $R$ to satisfy
\[ |R|=O(|x|^{-q})
\]
for some $q>n$. Under these assumptions
the ADM mass is well defined by the formula (see \cite{sc} for the $n$ dimensional case)
\[ m=\frac{1}{4(n-1)\o_{n-1}}\lim_{\s\to\infty}\int_{S_\s}\sum_{i,j}(g_{ij.i}-g_{ii,j})\nu_j\ d\xi(\s)
\]
where $S_\s$ is the euclidean sphere in the $x$ coordinates, $\o_{n-1}=Vol(S^{n-1}(1))$, and 
the unit normal and volume integral are with respect to the euclidean metric. The positive
mass theorem is as follows.
\begin{thm}
Assume that $M$ is an asymptotically flat manifold with $R\geq 0$. We then have that the ADM mass is nonnegative. Furthermore, if the mass is zero, then $M$ is isometric to $\R^n$.
\end{thm}

It is shown in Section 5 using results of \cite{sy3} to simplify the asymptotic behavior and
an observation of J. Lohkamp which allows us to compactify the manifold keeping the
scalar curvature positive. The result which is needed for compact manifolds follows.
\begin{thm} If $M_1$ is any closed manifold of dimension $n$, then $M_1\#T^n$ does not
have a metric of positive scalar curvature.
\end{thm}
Both of these theorems were known if either $n\leq 8$ or for any $n$ assuming the manifold
is a spin manifold. Actually for $n=8$ there may be isolated singularities, but in this dimension
a result of N. Smale \cite{sm} shows that there is a dense set of ambient metrics for which
the singularities do not occur. Using this result the  eight dimensional case can also be done
without dealing with singularities. In this paper we remove the dimensional and spin assumptions.

Finally we prove the following more precise theorem about compact manifolds with positive
scalar curvature.
\begin{thm}
Assume that $M$ is a compact oriented $n$-manifold with a metric of positive
scalar curvature. If $\a_1,\ldots,\a_{n-2}$ are classes in $H^1(M,\Z)$ with the property that the
class $\s_2$ given by
$\s_2=\a_{n-2}\cap\a_{n-3}\cap\ldots\a_1\cap[M]\in H_2(M,\Z)$ is nonzero, then the class 
$\s_2$ can be represented by a sum of smooth two spheres.
If $\a_{n-1}$ is any class in $H^1(M,\Z)$, then we must have $\a_{n-1}\cap\s_2=0$. In particular,
if $M$ has classes $\a_1,\ldots,\a_{n-1}$ with $\a_{n-1}\cap\ldots\cap\a_1\cap[M]\neq 0$,
then $M$ cannot carry a metric of positive scalar curvature.
\end{thm}

We also point out the recent series of papers by J. Lohkmp \cite{lo1}, \cite{lo2}, \cite{lo3},
and \cite{lo4}. These papers also present an approach to the high dimensional positive mass
theorem by extending the minimal hypersurface approach to all dimensions. Our approach
seems quite different both conceptually and technically, and is more in the classical spirit of
the calculus of variations. In any case we feel that, for such a fundamental result, it is of
value to have multiple approaches.

\bigskip

\section{\bf  Terminology and statements of main theorems}

We begin by introducing the notation involved in the construction of a {\it minimal $k$-slicing}; that is,
a nested family of hypersurfaces beginning with a smooth manifold $\Sig_n$ of dimension
$n$ and going down to $\Sig_k$ of dimension $k\leq n-1$. This consists of $\Sig_k\subset\Sig_{k+1}
\subset\ldots \subset \Sig_n$ where each $\Sig_j$ will be constructed as a volume minimizer
of a certain weighted volume in $\Sig_{j+1}$.
 
Let $\Sig_n$ be a properly embedded $n$-dimensional submanifold in an open set $\O$
contained in $\R^N$. We will consider a minimal slicing of $\Sig_n$ defined in an inductive manner.
First, let $u_n=1$, and let $\Sig_{n-1}$ be a volume minimizing hypersurface in $\Sig_n$. 
Of course, it may happen that
$\Sig_{n-1}$ has a singular set ${\S}_{n-1}$ which is a closed subset of Hausdorff dimension
at most $n-8$. On $\Sig_{n-1}$ we will construct a positive definite quadratic form $Q_{n-1}$
on functions by suitably modifying the index form associated to the second variation of
volume. We will then construct a positive function $u_{n-1}$ on $\Sig_{n-1}$ which is a 
least eigenfunction of $Q_{n-1}$. We then define $\rho_{n-1}=u_{n-1}u_n$, and we let $\Sig_{n-2}$
be a hypersurface in $\Sig_{n-1}$ which is a minimizer of the $\rho_{n-1}$-weighted volume
 $V_{\rho_{n-1}}(\Sig)=\int_\Sig\rho_{n-1}d\mu_{n-2}$ for an $n-2$ dimensional submanifold of
 $\Sig_{n-1}$ and we denote $\mu_j$ to be the Hausdorff $j$-dimensional measure. Inductively,
 assume that we have constructed a slicing down to dimension $k+1$; that is, we have a nested
 family of hypersurfaces, quadratic forms, and positive functions $(\Sig_j,Q_j,u_j)$ for 
 $j=k+1,\ldots,n$ such that $\Sig_j$ minimizes the $\rho_{j+1}$-weighted volume where 
 $\rho_{j+1}=u_{j+1}u_{j+2}\ldots u_n$, $Q_j$ is a positive definite quadratic form related to
 the second variation of the $\rho_{j+1}$-weighted volume (see (\ref{eqn:qform}) below), and $u_j$ is a lowest eigenfunction of $Q_j$ with eigenvalue $\lambda_j\geq 0$. We will always take $\lambda_j$ to be the lowest Dirichlet eigenvalue (if $\partial\Sig_j\neq 0$) of $Q_j$ with respect to the weighted $L^2$ norm and we take $u_j$ to be a corresponding eigenfunction. We will show in Section 3 that such $\lambda_j$ and $u_j$ exist. We then inductively
 construct $(\Sig_k,Q_k,u_k)$ by letting $\Sig_k$ be a minimizer of the $\rho_{k+1}$ weighted volume
 where $\rho_{k+1}=u_{k+1}u_{k+2}\ldots u_n$, $Q_k$ a positive definite quadratic form described
 below, and $u_k$ a positive eigenfunction of $Q_k$.
 
 Note that if $\Sig_j$ is a leaf in a minimal $k$-slicing, then choosing a unit
normal vector $\nu_j$ to $\Sig_j$ in $\Sig_{j+1}$ gives us an orthonormal basis 
$\nu_k,\nu_{k+1},\ldots,\nu_{n-1}$ for the normal bundle of $\Sig_k$ defined on the regular
set ${\cal R}_k$. Thus the second fundamental form of $\Sig_k$ in $\Sig_n$ consists of the
scalar forms $A_k^{\nu_j}=\langle A_k,\nu_j\rangle$for $j=k,\ldots,n-1$ and we have
$|A_k|^2=\sum_{j=k}^{n-1}|A_k^{\nu_j}|^2$.  

 Now if we have a minimal $k$-slicing, we let $g_k$ denote the metric induced on $\Sig_k$ from
$\Sig_n$, and we let $\hat{g}_k$ denote the metric $\hat{g}_k=g_k+\sum_{p=k}^{n-1}u_p^2dt_p^2$
on $\Sig_k\times(S^1)^{n-k}$ where we use $S^1$ to denote a circle of length $1$, and we denote
by $t_p$ a coordinate on the $p$th factor of $S^1$. We then note that the volume measure
of the metric $\hat{g}_k$ is given by $\rho_k d\mu_k$ where we have suppressed the $t_p$ variables since we will consider only objects which do not depend on them; for example, the $\rho_k$-weighted volume of $\Sig_k$ is the volume of the $n$-dimensional manifold $\Sig_k\times T^{n-k}$. We will need to introduce another metric $\tilde{g}_k$ on $\Sig_k\times(S^1)^{n-k-1}$. This is defined by 
$\tilde{g}_k= g_k+\sum_{p=k+1}^{n-1}u_p^2\ dt_p^2$. Note that $\tilde{g}_k$ is the metric 
induced on $\Sig_k\times(S^1)^{n-k-1}$ by $\hat{g}_{k+1}$. We also let $\tilde{A}_k$ denote the 
second fundamental form of $\Sig_k\times(S^1)^{n-k-1}$ in $(\Sig_{k+1}\times(S^1)^{n-k-1}, \hat{g}_{k+1})$. The following lemma computes this second fundamental form.
\begin{lem} \label{lem:2ff} We have $\tilde{A}_k=A_k^{\nu_k}-\sum_{p=k+1}^{n-1}u_p\nu_k(u_p)dt_p^2$,
and the square length with respect to $\tilde{g}_k$ is given by 
$|\tilde{A}_k|^2=|A_k^{\nu_k}|^2+\sum_{p=k+1}^{n-1}(\nu_k(log\ u_p))^2$.
\end{lem}
\begin{pf} If we consider a hypersurface $\Sig$ in a Reimannian manifold with unit normal
$\nu$, then we can consider the parallel hypersurfaces parametrized on $\Sig$ by
$F_\e(x)=\exp(\e\nu(x))$ for small $\e$ and $x\in\Sig$. We then have a family of induced metrics $g_\e$ from $F_\e$  on $\Sig$, and the second fundamental form is given by $A=-\frac{1}{2}\dot{g}$ 
where $\dot{g}$ denotes the $\e$ derivative of $g_\e$ at $\e=0$.

If we let $\exp$ denote the exponential map of $\Sig_k$ in $\Sig_{k+1}$, then since $\Sig_{k+1}$
is totally geodesic in $\Sig_{k+1}\times T^{n-k-1}$, we have 
\[ F_\e(x,t)=(\exp(\e\nu_k(x),t) 
\]
for $(x,t)\in \Sig_k\times T^{n-k-1}$, and the induced family of metrics is given by 
\[ \tilde{g}_\e=(g_k)_\e+\sum_{p=k+1}^{n-1}(u_p(\exp(\e\nu_k))^2\ dt_p^2.
\] 
Thus we have 
\[ \dot{\tilde{g}}=-2A_k^{\nu_k}+2\sum_{p=k+1}^{n-1}u_p\nu_k(u_p)\ dt_p^2
\]
since $A_k^{\nu_k}$ is the second fundamental form of $\Sig_k$ in $\Sig_{k+1}$. It follows that
$\tilde{A}_k=A_k^{\nu_k}-\sum_{p=k+1}^{n-1}u_p\nu_k(u_p)dt_p^2$, and taking the square norm
with respect to the metric $\tilde{g}_k$ then gives the desired formula for $|\tilde{A}_k|^2$.
\end{pf}
 
We now describe the choice we will make for $Q_j$. Let $S_j$ be the second variation
form for the weighted volume $V_{\rho_{j+1}}$ at $\Sig_j$, and define
\begin{eqnarray}
\label{eqn:qform}
Q_j(\varphi,\varphi)&=&S_j(\varphi,\varphi)+\frac{3}{8}\int_{\Sig_j}
(|\tilde{A}_j|^2\nonumber\\ &+&\frac{1}{3n} \sum_{p=j+1}^n(|\nabla_jlog\ u_p|^2+|\tilde{A}_p|^2))\varphi^2
\rho_{j+1}\ d\mu_j
\end{eqnarray}
where, for now, $\varphi$ is a function supported in the regular set ${\cal R}_j$ and we define 
$\tilde{A}_n=0,\ u_n=1$. We will discuss an extended domain for $Q_j$ in the Section 3.

Up to this point our discussion is formal because we have not discussed issues related to
 the singularities of the $\Sig_j$ in a minimal slicing. We first define the {\it regular set}, ${\cal R}_j$
 of $\Sig_j$ to be the set of points $x$ for which there is a neighborhood of $x$ in $\R^N$ in which
 all of $\Sig_j,\Sig_{j+1},\ldots \Sig_n$ are smooth embedded submanifolds of $\R^N$. The {\it singular
 set}, ${\S}_j$ is then defined to be the complement of ${\cal R}_j$ in $\Sig_j$. Thus ${\S}_j$ is
 a closed set by definition. The following result follows from the standard minimizing hypersurface
 regularity theory. In this paper $dim(A)$ always refers to the Hausdorff dimension of a subset
 $A\subset \R^N$.
 \begin{prop}
\label{prop:topreg}
For $j\leq n-1$ we have $dim({\S}_j\sim {\S}_{j+1})\leq j-7$, and in particular we 
have $dim({\S}_{n-1})\leq n-8$. 
\end{prop}
In light of this result, we see that our main task in controlling singularities is to control
the size of the set ${\S}_j\cap {\S}_{j+1}$. We will do this by extending the minimal
hypersurface regularity theory to this slicing setting. In order to do this we need to establish
the relevant compactness and tangent cone properties and this requires establishing
suitable bounds on the slicings. To begin this process we make the following definition.
\begin{defn} For a constant $\Lambda>0$, a {\bf $\Lambda$-bounded minimal $k$-slicing} 
is a minimal $k$-slicing satisfying the following bounds
$$ \lambda_j\leq \Lambda,\ Vol_{\rho_{j+1}}(\Sig_j)\leq \Lambda,\ \int_{\Sig_j}(1+|A_j|^2+
\sum_{p=j+1}^n|\nabla_jlog\ u_p|^2)u_j^2\rho_{j+1}\ d\mu_j\leq \Lambda
$$ 
for $j=k,k+1,\ldots n-1$, where $\mu_j$ is Hausdorff measure, $\nabla_j$ is taken on (the 
regular set of) $\Sig_j$, and $A_j$ is the second fundamental form of $\Sig_j$ in $\R^N$.
\end{defn}
The minimal $k$-slicings we will consider in this paper will always be $\Lambda$-bounded
for some $\Lambda$. We have the following regularity theorem.

\begin{thm}
\label{thm:reg} Given any $\Lambda$-bounded minimal $k$-slicing, we have for each
$j=k,k+1,\ldots, n-1$ the bound on the singular set $dim({\S}_j)\leq j-3$. 
\end{thm}

We now formulate an existence theorem for minimal $k$-slicings in $\Sig_n$. We consider the
case in which $\Sig_n$ is a closed oriented manifold. We assume that there is closed
oriented $k$-dimensional manifold $X^k$ and a smooth map $F:\Sig_n\to X\times T^{n-k}$ of
non-zero degree $s$. We let $\Omega$ denote a $k$-form of $X$ with $\int_X\Omega=1$, and we denote by $dt^{k+1},\ldots dt^n$ the basic one forms on $T^{n-k}$ where we assume the periods are
equal to one. We introduce the notation $\Theta=F^*\Omega$ and $\omega^p=F^*(dt^p)$
for $p=k+1,\ldots, n$.

We can now state our first existence theorem. A more refined existence theorem is given
by Theorem \ref{thm:exst2} which we will not state here.
\begin{thm}
\label{thm:exst} For a manifold $M=\Sig_n$ as described above, there is a $\Lambda$-bounded,
partially regular, minimal $k$-slicing Moreover, if $k\leq j\leq n-1$ and $\Sig_j$ is regular, then 
$\int_{\Sig_j}\Theta\wedge\omega^{k+1}\wedge\ldots\wedge\omega^j=s$.
\end{thm}

The proofs of Theorems \ref{thm:reg} and \ref{thm:exst} will be given in Sections 3 and 4. 
In the
remainder of this section we discuss the quadratic forms $Q_j$ in more detail and derive 
important geometric consequences for minimal $1$-slicings and $2$-slicings under the 
assumption that $\Sig_n$ has positive scalar curvature. Consequences of these results, which
are the main geometric theorems of the paper, will be given in Section 5. 

Recall that in general if $\Sig$ is a stable two-sided (trivial normal bundle) minimal hypersurface in a Riemannian manifold $M$, then we may choose a globally defined unit normal vector $\nu$,
and we may parametrize normal deformations by functions $\varphi\cdot\nu$. The second variation
of volume then becomes the quadratic form
\begin{equation}
\label{eqn:secvar} 
S(\varphi,\varphi)=\int_\Sig[|\nabla\varphi|^2-\frac{1}{2}(R_M-R_\Sig+|A|^2)\varphi^2]\ d\mu
\end{equation}
where $R_M$ and $R_\Sig$ are the scalar curvature functions of $M$ and $\Sig$ and $A$
denotes the second fundamental form of $\Sig$ in $M$.

We have the following result which
computes the scalar curvature $\tilde{R}_k$ of $\tilde{g}_k$.
\begin{lem} 
\label{lem:rcalc} The scalar curvature of the metric $\tilde{g}_k$ is given by
\[ \tilde{R}_k=R_k-2\sum_{p=k+1}^{n-1}u_p^{-1}\Delta_ku_p-2\sum_{k+1\leq p<q\leq n-1}
\langle\nabla_k log\ u_p,\nabla_k log\ u_q\rangle
\]
where $\Delta_k$ and $\nabla_k$ denote the Laplace and gradient operators with respect to
$g_k$.
\end{lem}
\begin{pf} The calculation is a finite induction using the formula
\[ \tilde{R}=R-2u^{-1}\Delta u
\]
for the scalar curvature of the metric $\tilde{g}=g+u^2dt^2$. 

For $j=k,\ldots,n-1$ 
Let $\bar{g}_j=g_k+\sum_{p=j}^{n-1}u_p^2dt_p^2$. Note that $\bar{g}_k=\hat{g}_k$
and $\bar{g}_{k+1}=\tilde{g}_k$. We prove the formula
\[ \bar{R}_j=R_k-2\sum_{p=j}^{n-1}u_p^{-1}\Delta_k u_p-2\sum_{j\leq p<q\leq n-1}
\langle\nabla_k log\ u_p,\nabla_k log\ u_q\rangle
\]
by a finite reverse induction on $j$. First note that for $j=n-1$ the formula follows from
the one above. Now assume the formula is correct for $\bar{g}_{j+1}$
We then apply the formula above to obtain
\[ \bar{R}_j=\bar{R}_{j+1}-2u_j^{-1}\bar{\Delta}_j u_j. 
\]
Since $u_j$ does not depend on the extra variables $t_p$, we have 
\[ u_j^{-1}\bar{\Delta}_j u_j=u_j^{-1}\rho_j^{-1}div_k(\rho_j\nabla_k u_j)=
u_j^{-1}\Delta_k u_j
+\sum_{p=j+1}^{n-1}\langle\nabla_k log\ u_p, \nabla_k log\ u_j\rangle
\]
where as above $\rho_j=u_{j+1}\cdots u_{n-1}$. The statement now follows from
the inductive assumption. Since $\bar{g}_{k+1}=\tilde{g}_k$, we have proven the required
statement.
\end{pf}
We now consider consequences of having a minimal $k$-slicing of a manifold of positive
scalar curvature. 
\begin{thm}
\label{thm:eval} Assume that the scalar curvature of $\Sig_n$ is bounded below by a constant
$\kappa$. If $\Sig_k$ is a leaf in a minimal $k$-slicing, then we have the following scalar curvature 
formula and eigenvalue estimate 
\[  \hat{R}_k= R_n+2\sum_{p=k}^{n-1}\lambda_p+\frac{1}{4}\sum_{p=k}^{n-1}(|\tilde{A}_p|^2
-\frac{1}{n}\sum_{q=p+1}^n(|\nabla_plog\ u_q|^2+|\tilde{A}_q|^2))
\] 
\[ \int_{\Sig_k}(\kappa+\frac{3}{4}\sum_{j=k+1}^n|\nabla_klog\ u_j|^2-R_k)\varphi^2\ d\mu_k\leq 4\int_{\Sig_k}|\nabla_k\varphi|^2\ d\mu_k
\]
where $\varphi$ is any smooth function with compact support in ${\cal R}_k$.
\end{thm}
\begin{pf} First note that from (\ref{eqn:qform}) and (\ref{eqn:secvar}) we have
\begin{eqnarray*}
Q_j(\varphi,\varphi)&=&\int_{\Sig_j}[|\nabla_j\varphi|^2-
\frac{1}{2}(\hat{R}_{j+1}-\tilde{R}_j)\varphi^2 \\
&-&\frac{1}{8}(|\tilde{A}_j|^2-\frac{1}{n}\sum_{p=j+1}^n(|\nabla_jlog\ u_p|^2+|\tilde{A}_p|^2))\varphi^2]\rho_{j+1}\ d\mu_j,
\end{eqnarray*}
and therefore $u_j$ satisfies the equation $L_ju_j=-\lambda_ju_j$ where
\begin{equation} 
\label{eqn:operator} L_j=\tilde{\Delta}_j+\frac{1}{2}(\hat{R}_{j+1}-\tilde{R}_j)+\frac{1}{8}(|\tilde{A}_j|^2-\frac{1}{n}\sum_{p=j+1}^n(|\nabla_jlog\ u_p|^2+|\tilde{A}_p|^2)).
\end{equation}

We derive the scalar curvature formula by a finite downward induction beginning with $k=n-1$. 
In this case 
the eigenvalue estimates follow from the standard stability inequality (\ref{eqn:secvar}) since 
$\rho_n=u_n=1$ and $\tilde{R}_{n-1}=R_{n-1}$. We also have from Lemma \ref{lem:rcalc} that 
$\hat{R}_{n-1}=R_{n-1}-2u_{n-1}^{-1}\Delta_{n-1}u_{n-1}$. The equation satisfied by $u_{n-1}$ is 
\[ \Delta_{n-1}u_{n-1}+\frac{1}{2}(R_n-R_{n-1})u_{n-1}+\frac{1}{8}|\tilde{A}_{n-1}|^2u_{n-1}=-\lambda_{n-1}u_{n-1}
\] 
and so we have 
$\hat{R}_{n-1}=R_n+2\lambda_{n-1}+\frac{1}{4}|\tilde{A}_{n-1}|^2$. This proves the result for $k=n-1$.

Now we assume the conclusions are true for integers $k$ and larger, and we will derive them 
for $k-1$. We first observe that $\hat{g}_{k-1}=\tilde{g}_{k-1}+u_{k-1}^2\ dt_{k-1}^2$ and so
$\hat{R}_{k-1}=\tilde{R}_{k-1}-2u_{k-1}^{-1}\tilde{\Delta}_{k-1}u_{k-1}$. On the other hand
from (\ref{eqn:operator}) applied with $j=k-1$ we see that $u_{k-1}$ satisfies the equation
\begin{eqnarray*}
\tilde{\Delta}_{k-1}u_{k-1}&+&\frac{1}{2}(\hat{R}_k-\tilde{R}_{k-1})u_{k-1}+\frac{1}{8}(|\tilde{A}_{k-1}|^2
\\ &-&\frac{1}{n}\sum_{p=k}^n(|\nabla_{k-1}log\ u_p|^2+|\tilde{A}_p|^2))u_{k-1}=-\lambda_{k-1}u_{k-1}.
\end{eqnarray*}
Substituting this above we have 
\begin{eqnarray*}
 \hat{R}_{k-1}&=&\tilde{R}_{k-1}+2[\lambda_{k-1}+\frac{1}{2}(\hat{R}_k-\tilde{R}_{k-1}) \\
&+&\frac{1}{8}(|\tilde{A}_{k-1}|^2-\frac{1}{n}\sum_{q=k}^n(|\nabla_{k-1}log\ u_q|^2+|\tilde{A}_q|^2))],
\end{eqnarray*}
so we have
\[ \hat{R}_{k-1}=2\lambda_{k-1}+\hat{R}_k+\frac{1}{4}(|\tilde{A}_{k-1}|^2-\frac{1}{n}\sum_{q=k}^n(|\nabla_{k-1}log\ u_q|^2+|\tilde{A}_q|^2)).
\]
Using the inductive hypothesis  we get the desired formula 
\[\hat{R}_{k-1}= R_n+2\sum_{p=k-1}^{n-1}\lambda_p+\frac{1}{4}\sum_{p=k-1}^{n-1}(|\tilde{A}_p|^2
-\frac{1}{n}\sum_{q=p+1}^n(|\nabla_p log\ u_q|^2+|\tilde{A}_q|^2)).
\]

Now observe that 
\begin{eqnarray*}
 \sum_{p=k}^{n-1}(n|\tilde{A}_p|^2&-&\sum_{q=p+1}^n(|\nabla_p log\ u_q|^2+|\tilde{A}_q|^2)) \\
&\geq& \sum_{p=k}^{n-1}(\sum_{r=k}^n|\tilde{A}_r|^2-\sum_{q=p+1}^n(|\nabla_p log\ u_q|^2+|\tilde{A}_q|^2)) \\
&\geq& \sum_{p=k}^{n-1}\sum_{q=p+1}^n(\sum_{r=k}^{p-1}(\nu_rlog\ (u_q))^2-|\nabla_plog\ u_q|^2)\\
&=&-\sum_{p=k}^{n-1}\sum_{q=p+1}^n|\nabla_{k-1}log\ u_q|^2\geq -n\sum_{q=k}^n|\nabla_{k-1}log\ u_q|^2.
\end{eqnarray*} 
This formula implies that for each $k$ we have 
\begin{equation}\label{eqn:scbound} \hat{R}_k\geq\kappa-1/4\sum_{j=k}^n|\nabla_{k-1}log\ u_j|^2
\end{equation}
and so the following eigenvalue estimate follows from (\ref{eqn:secvar}) 
\[ \int_{\Sig_k}(\kappa-\frac{1}{4}\sum_{j=k+1}^n|\nabla_k log\ u_j|^2-\tilde{R}_k)\varphi^2\rho_{k+1}\ d\mu_k
\leq 2\int_{\Sig_k}|\nabla_k\varphi|^2\rho_{k+1}\ d\mu_k
\]
The remainder of the proof derives the eigenvalue estimate from this one. Since $\varphi$
is arbitrary we may replace $\varphi$ by $\varphi(\rho_{k+1})^{1/2}$ to obtain
\begin{eqnarray*}
\int_{\Sig_k}(\kappa-\frac{1}{4}\sum_{j=k+1}^n|\nabla_klog\  u_j|^2-\tilde{R}_k)\varphi^2\ d\mu_k&\leq&
2\int_{\Sig_k}|\nabla_k(\varphi/\sqrt{\rho_{k+1}})|^2\rho_{k+1}\ d\mu_k \\
&\leq& 4\int_{\Sig_k}|\nabla_k(\varphi/\sqrt{\rho_{k+1}})|^2\rho_{k+1}\ d\mu_k
\end{eqnarray*} 
where we used the inequality $2\leq 4$. After expanding, the term on the right becomes
\[ 4\int_{\Sig_k}(|\nabla_k\varphi|^2-\varphi\langle\nabla_k\varphi,\nabla_klog\ \rho_{k+1}\rangle
+1/4\varphi^2|\nabla_k log\ \rho_{k+1}|^2)\ d\mu_k.
\]
Rewriting the middle term in terms of $\nabla_k(\varphi)^2$ and integrating by parts the term 
becomes
\[ 4\int_{\Sig_k}(|\nabla_k\varphi|^2+1/2\varphi^2[\sum_{p=k+1}^{n-1}(u_p^{-1}\Delta_k u_p
-|\nabla_k log\ u_p|^2)+1/2|\nabla_k log\ \rho_{k+1}|^2])\ d\mu_k.
\]
Now recall from Lemma \ref{lem:rcalc} that
\[ \tilde{R}_k=R_k-2\sum_{p=k+1}^{n-1}u_p^{-1}\Delta_k u_p-2\sum_{k+1\leq p<q\leq n-1}
\langle\nabla_k log\ u_p,\nabla_k log\ u_q\rangle.
\]
Thus we see that the terms involving $\Delta_k u_p$ cancel out, and note also that
\[ |\nabla_k log\ \rho_{k+1}|^2=\sum_{p=k+1}^{n-1}|\nabla_k\ log\ u_p|^2+2\sum_{k+1\leq p<q\leq n-1}
\langle\nabla_k log\ u_p,\nabla_k log\ u_q\rangle
\]
so the second term also cancels. Thus we are left with
\begin{eqnarray*}
\int_{\Sig_k}(\kappa-\frac{1}{4}\sum_{j=k+1}^n|\nabla_klog\ u_j|^2&-&R_k)\varphi^2\ d\mu_k \\
&\leq& 4\int_{\Sig_k}(|\nabla_k\varphi|^2-\frac{1}{4} \sum_{j=k+1}^n|\nabla_k\ log\ u_j|^2)\ d\mu_k.
\end{eqnarray*}
This gives the desired eigenvalue estimate.
\end{pf}
This theorem will be central to the regularity proof in the next section and it also has an
important geometric consequence which is the main tool in the applications of Section 5.
\begin{thm}
\label{thm:12slicing} Assume that $R_n\geq \kappa>0$. If $\Sig_k$ is regular, then
$(\Sig_k,g_k)$ is a Yamabe positive conformal manifold. If $\Sig_2$ lies in a minimal 
$2$-slicing, $\Sig_2$ is regular, and $\partial\Sig_2=0$, then each connected
component of $\Sig_2$ is homeomorphic to the two sphere. If $\Sig_1$ lies in a minimal
$1$-slicing and $\Sig_1$ is regular, then each component of $\Sig_1$ is an arc of length
at most $2\pi/\sqrt{\kappa}$.
\end{thm}
\begin{pf} Recall that the condition that $g_k$ be Yamabe positive is that the lowest eigenvalue
of the conformal Laplacian $-\Delta_k+c(k)R_k$ be positive where $c(k)=\frac{k-2}{4(k-1)}$. In 
variational form this condition says
\[ -\int_{\Sig_k}R_k\varphi^2\ d\mu_k<c(k)^{-1}\int_{\Sig_k}|\nabla_k\varphi|^2\ d\mu_k
\]
for all nonzero functions $\varphi$ which vanish on $\partial\Sig_k$ (if $\Sig_k$ has a boundary).
Since $4<c(k)^{-1}$ we see that this follows from the eigenvalue estimate of Theorem \ref{thm:eval}.

Now consider $\Sig_2$, and apply the eigenvalue estimate of Theorem
\ref{thm:eval} with $\varphi=1$ to a component $S$ of $\Sig_2$ to see that 
$\int_S R_2\ d\mu_2>0$. It then follows from the Gauss-Bonnet Theorem that $S$ is
homeomorphic to the two sphere (note that $S$ is orientable).

Finally, it $\gamma$ is a connected component of $\Sig_1$ of length $l$, then the
eigenvalue estimate of Theorem \ref{thm:eval} implies that the lowest Dirichlet eigenvalue
of $\gamma$ is at least $\kappa/4$. Thus $\kappa/4\leq \pi^2/l^2$ and $l\leq 2\pi/\sqrt{\kappa}$
as claimed.
\end{pf}
\bigskip

\section{\bf Compactness and regularity of minimal $k$-slicings}
The main goal of this section is to prove Theorem \ref{thm:reg}. In order
to do this we first must clarify some analytic issues concerning the domain of the quadratic form
$Q_j$. We let $L^2(\Sig_j)$ denote the space of square integrable 
functions on $\Sig_j$ with respect to the measure $\rho_{j+1}\mu_j$. We let 
\[ \|\varphi\|^2_{0,j}=\int_{\Sig_j}\varphi^2 \rho_{j+1}\ d\mu_j
\]
denote the square norm on $L^2_{\Sig_j}$. We introduce some notation, defining $P_j$ to be
the function defined on $\Sig_j$
\[ P_j=|A_j|^2+\sum_{p=j+1}^n|\nabla_jlog\ u_p|^2.
\]

We will say that a minimal $k$-slicing in an open set $\Omega$ is {\it partially regular} if 
$dim({\S}_j)\leq j-3$ for $j=k,\ldots,n-1$. It follows from Proposition \ref{prop:topreg} that
if the $(k+1)$-slicing associated to a minimal $k$-slicing is partially regular, then
$dim({\S}_k)\leq min\{dim({\S}_{k+1}),k-7\}\leq k-2$.

For functions $\varphi$ which are Lipschitz (with respect to ambient distance) on $\Sig_j$ with compact support in ${\cal R}_j\cap\bar{\O}$, we define a square norm by
\[ \|\varphi\|_{1,j}^2=\|\varphi\|^2_{0,j}+ \int_{\Sig_j}(|\nabla_j\varphi|^2+P_j\varphi^2)\rho_{j+1}\ d\mu_j.
\] 
We let ${\cal H}_j$ denote the Hilbert space which is the completion with respect to this norm. Note
that functions in ${\cal H}_j$ are clearly locally in $W_{1,2}$ on ${\cal R}_j$. We will
assume from now on that $u_j\in {\cal H}_j$ for $j\geq k$; in fact, we take this as part of the 
definition of a bounded minimal $k$-slicing. We define ${\cal H}_{j,0}$ to be the closed subspace of ${\cal H}_j$ consisting of the completion of the Lipschitz functions with compact support in ${\cal R}_j\cap\Omega$.
In order to handle boundary effects we also assume that there is a larger domain $\O_1$ which contains $\bar{\O}$ as a compact subset and that the $k$-slicing is defined and boundaryless in $\O_1$. Note that this is automatic if $\partial\Sig_j=\phi$. Thus ${\cal H}_{j,0}$ consists of those functions in ${\cal H}_j$ with $0$ boundary data on $\Sig_j\cap\partial\O$. The existence of eigenfunctions $u_j$ in this 
space will be discussed in the next section. The following estimate of the $L^2(\Sig_k)$ norm
near the singular set will be used both in this section and the next. The result may be thought
of as a non-concentration result for the weighted $L^2$ norm near the singular set in case
the ${\cal H}_k$ norm is bounded.
\begin{prop}
\label{prop:l2con} Let ${\cal S}$ be a closed subset of $\O_1$ with zero $(k-1)$-dimensional Hausdorff measure. Let $\Sig_k$ be a member of a bounded minimal $k$-slicing such that $\Sig_{k+1}$ is partially regular in $\O_1$. For any $\eta>0$ there
exists an open set $V\subset \O_1$ containing ${\cal S}\cap\bar{\O}$ such that whenever
${\cal S}_k\cap\bar{\O}\subset V$ we have the following 
estimate
\[ \int_{\Sig_k\cap V}\varphi^2\rho_{k+1}\ d\mu_k\leq
\eta\int_{\Sig_k\cap\O}[|\nabla_k\varphi|^2+(1+P_k)\varphi^2]\rho_{k+1}\ d\mu_k
\]
for all $\varphi\in {\cal H}_{k,0}$. 
\end{prop}
\begin{pf} Let $\e>0,\ \d>0$ be given. We may choose a finite covering of the compact set 
${\cal S}\cap \bar{\O}$ by balls $B_{r_\a}(x_\a)$ with $r_\a\leq \d/5$ such
\[ \sum_\a r_\a^{k-1}\leq \e.
\]
We let $V$ denote the union of the balls, $V=\cup_\a B_{r_\a}(x_\a)$. 

Assume that ${\cal S}_k\cap\bar{\O}\subset V$ and let $\varphi\in {\cal H}_{k,0}$.
We may extend $\varphi$ to $\Sig_k\cap\O_1$ be taking $\varphi=0$ in $\O_1\sim\O$.
 By a standard first variation argument for submanifolds of $\R^N$, for a nonnegative function we have
\begin{eqnarray*}
k\int_{\Sig_k\cap B_r}\varphi^2\rho_{k+1}\ d\mu_k&\leq& r\int_{\Sig_k\cap B_r}(|\nabla_k
\varphi^2\rho_{k+1}|+|H_k|\varphi^2\rho_{k+1})\ d\mu_k \\
&+&r\int_{\Sig_k\cap\partial B_r} \varphi^2\rho_{k+1}\ d\mu_{k-1}.
\end{eqnarray*}
Let $L_\a(r)=\int_{\Sig_k\cap B_r(x_\a)}\varphi^2\rho_{k+1}\ d\mu_k$ and
\[M_\a(r)=\int_{\Sig_k\cap B_r(x_\a)}(|\nabla_k
(\varphi^2\rho_{k+1})|+|H_k|\varphi^2\rho_{k+1})\ d\mu_k.
\]
The above inequality then implies
\[ kL_\a(r)\leq rM_\a(r)+r\frac{d}{dr}(L_\a(r)).
\]
Now for any $\a$ and a small constant $\e_0$ we consider two cases: (1) There exists $r$ with 
$r_\a\leq r\leq \d/5$ such that the inequality
\[ \e_0L_\a(5r)\leq rM_\a(r).
\]
We denote such a choice of $r$ by $r_\a'$. Secondly, we have case (2) For all $r$ with 
$r_\a\leq r\leq \d/5$ we have
\[ rM_\a(r)< \e_0L_\a(5r).
\]
The collection of $\a$ for which the first case holds will be labeled $A_1$, and that for which the 
second holds $A_2$. We will handle the two cases separately.

For the collection of balls with radius $r_\a'$ indexed by $A_1$ we may apply the five times 
covering lemma to extract a subset $A_1'\subseteq A_1$ for which the balls in $A_1'$ are 
disjoint and such that
\[ V_1\equiv \cup_{\a\in A_1}B_{r_\a}(x_\a)\subseteq  \cup_{\a\in A_1}B_{r_\a'}(x_\a)\subseteq 
\cup_{\a\in A_1'}B_{5r_\a'}(x_\a).
\]
From the inequality of case (1) above applied for $\a\in A_2'$ we have
\[ L_\a(r_\a)\leq L_\a(5r_\a')\leq \e_0^{-1}r_\a'M_\a(r_\a')\leq \e_0^{-1}\d M_\a(r_\a').
\] 
Summing over $\a\in A_1$ and using disjointness of the balls we have
\begin{equation}
\label{eqn:case1} \int_{\Sig_k\cap V_1}\varphi^2\rho_{k+1}\ d\mu_k\leq 
\e_0^{-1}\d\int_{\Sig_k\cap\O}(|\nabla_k \varphi^2\rho_{k+1}|+|H_k|\varphi^2\rho_{k+1})\ d\mu_k.
\end{equation}

Now for $\a\in A_2$ we have 
\[ kL_\a(r)\leq \e_0L_\a(5r)+r\frac{d}{dr}(L_\a(r)) 
\]
for $r_\a\leq r\leq \d/5$. For $j=0,1,2,\ldots$ define $\s_j=5^jr_\a$ and let $p$ be the
positive integer such that $\s_{p-1}<\d/5\leq \s_p$. We define $\Lambda_j$ by 
$\Lambda_j=L_\a(\s_j)$ for $j=0,1,\ldots,p$. For $\s_j\leq r\leq \s_{j+1}$ we then have
\[ kL_\a(r)\leq \e_0\Lambda_{j+2}\Lambda_j^{-1}L_\a(r)+r\frac{d}{dr}(L_\a(r)). 
\]
Integrating we find
\[ \Lambda_{j+1}\Lambda_j^{-1}\geq 5^{k-\e_0\Lambda_{j+2}\Lambda_j^{-1}}.
\]
Setting $R_j=\Lambda_{j+1}\Lambda_j^{-1}$ we have shown
\[ R_j\geq 5^{k-\e_0R_jR_{j+1}}.
\]  
Now if $R_j\leq 5^{k-1}$ then we would have $5^{k-1}\geq 5^{k-\e_0R_jR_{j+1}}$
which in turn implies $\e_05^{k-1}R_{j+1}\geq \e_0R_jR_{j+1}\geq 1$. Thus if we
choose $\e_0=5^{-3k+3}$ we find $R_{j+1}\geq 5^{2(k-1)}$ and hence it follows
that $R_jR_{j+1}\geq 5^{2(k-1)}$. Thus we have shown that for any $j=0,1,\ldots, p-1$
we either have $R_j\geq 5^{k-1}$ or $R_jR_{j+1}\geq 5^{2(k-1)}$. This implies that
$\Lambda_p\Lambda_0^{-1}\geq 5^{(p-1)(k-1)}\geq 5^{1-k}(\d/r_\a)^{k-1}$ and therefore
we have $L_\a(r_\a)\leq c(r_\a/\d)^{k-1}L_\a(\s_p)$ for each $\a\in A_2$. Summing this over
these $\a$ and using the choice of the covering we have
\[ \int_{\Sig_k\cap V_2}\varphi^2\rho_{k+1}\ d\mu_k\leq c\e\d^{1-k}\int_{\Sig_k\cap\O}\varphi^2\rho_{k+1}\ d\mu_k.
\]
Combining this with (\ref{eqn:case1}) we finally obtain
\[  \int_{\Sig_k\cap V}\varphi^2\rho_{k+1}\ d\mu_k\leq c\e\d^{1-k}\int_{\Sig_k\cap\O}\varphi^2\rho_{k+1}\ d\mu_k+c\d\int_{\Sig_k\cap\O}(|\nabla_k \varphi^2\rho_{k+1}|+|H_k|\varphi^2\rho_{k+1})\ d\mu_k.
\]
since we have now fixed $\e_0$. We can estimate the second term on the right using
\[ |\nabla_k \varphi^2\rho_{k+1}|+|H_k|\varphi^2\rho_{k+1}\leq (\varphi^2+|\nabla_k\varphi|^2)\rho_{k+1}
+\frac{1}{2}\varphi^2(2+|\nabla_k\log\ \rho_{k+1}|^2+|H_k|^2)\rho_{k+1}.
\]
This implies the bound
\[  \int_{\Sig_k\cap V}\varphi^2\rho_{k+1}\ d\mu_k\leq c(\e\d^{1-k}+\d)\int_{\Sig_k\cap\O}\varphi^2\rho_{k+1}\ d\mu_k+c\d\int_{\Sig_k\cap\O}[|\nabla_k \varphi|^2+P_k\varphi^2]\rho_{k+1}\ d\mu_k.
\]
The desired conclusion now follows by choosing $\d$ so that $c\d=\eta/2$ and then choosing
$\e$ so that $c\e\d^{1-k}=\eta$. This completes the proof.
\end{pf}

The following coercivity bound will be useful both in this section and in the next. We assume 
here that we have a partially regular minimal $k$-slicing. 
\begin{prop}
\label{prop:coercive} Assume that our $k$-slicing is bounded. There is a constant $c$ such that 
for $\varphi\in {\cal H}_{k,0}$ we have
\begin{equation*}
c^{-1}\int_{\Sig_k}[|\nabla_k\varphi|^2+(P_k+|\nabla_klog\ u_k|^2)\varphi^2] \rho_{k+1}\ d\mu_k
\leq Q_k(\varphi,\varphi)+\int_{\Sig_k}\varphi^2\rho_{k+1}\ d\mu_k.
\end{equation*}
Moreover we have the bound
\[ c^{-1}\int_{\Sig_k}(|\nabla_k(\varphi\sqrt{\rho_{k+1}})|^2+|A_k|^2\varphi^2\rho_{k+1})\ d\mu_k\leq 
Q_k(\varphi,\varphi)+\int_{\Sig_k}\varphi^2\rho_{k+1}\ d\mu_k.
\]
\end{prop}
\begin{pf} We can see from (\ref{eqn:qform}) that
\[ Q_k(\varphi,\varphi)\geq S_k(\varphi,\varphi)+\frac{1}{8n}\int_{\Sig_k}(\sum_{p=k}^n|\tilde{A}_p|^2+\sum_{p=k+1}^n|\nabla_k\log u_p|^2)\varphi^2\rho_{k+1}d\mu_k.
\]
Using the stability of $\Sig_k$ we have
\begin{equation} \label{eqn:qbound}
Q_k(\varphi,\varphi)\geq \frac{1}{8n}\int_{\Sig_k}(\sum_{p=k}^n|\tilde{A}_p|^2+\sum_{p=k+1}^n|\nabla_k\log u_p|^2)\varphi^2\rho_{k+1}d\mu_k.
\end{equation}
Finally we use Lemma \ref{lem:2ff} to conclude that (note that $\tilde{A}_n=0$)
\[ \sum_{p=k}^n|\tilde{A}_p|^2\geq \sum_{p=k}^{n-1}|A_p^{\nu_p}|^2\geq \sum_{p=k}^{n-1}|A_k^{\nu_p}|^2=|A_k|^2,
\]
and thus we have
\[ Q_k(\varphi,\varphi)\geq \frac{1}{8n}\int_{\Sig_k}P_k\varphi^2\rho_{k+1}\ d\mu_k.
\]
Recall that $S_k(\varphi,\varphi)=\int_{\Sig_k}(|\nabla_k\varphi|^2-q_k\varphi^2)\rho_{k+1}\ d\mu_k$ where
\[ q_k=\frac{1}{2}(|\tilde{A}_k|^2+\hat{R}_{k+1}-\tilde{R}_k)
\]
where $\hat{R}_{k+1}$ is given in Theorem \ref{thm:eval} and $\tilde{R}_k$ is given in Lemma
\ref{lem:rcalc}. We will need an upper bound on $q_k$, so we first see from Theorem \ref{thm:eval}
with $k$ replace by $k+1$
\[ q_k\leq c+\frac{1}{2}\sum_{p=k}^{n-1}|\tilde{A}_p|^2-\frac{1}{2}\tilde{R}_k
\]
where the constant bounds the curvature of $\Sig_n$ and the eigenvalues. Now from Lemma \ref{lem:rcalc} we can obtain the bound
\[ -\frac{1}{2}\tilde{R}_k\leq \frac{1}{2}|R_k|+\sum_{p=k+1}^{n-1}|\nabla_k\log u_p|^2+div_k({\cal X}_k)
\]
where ${\cal X}_k=\sum_{p=k+1}^{n-1}\nabla_k log\ u_p$. We observe that the Gauss equation
implies that $|R_k|\leq c(1+|A_k|^2)$, and so we have
\[ q_k\leq c+c\sum_{p=k}^{n-1}|\tilde{A}_p|^2+\sum_{p=k+1}^{n-1}|\nabla_k\log u_p|^2+div_k({\cal X}_k)
\]

Now observe that $Q_k\geq S_k$ and so we have
\[  \int_{\Sig_k}(|\nabla_k\varphi|^2+\frac{1}{8n}P_k\varphi^2)\rho_{k+1}\ d\mu_k\leq 2Q_k(\varphi,\varphi)
+\int_{\Sig_k}q_k\varphi^2\rho_{k+1}\ d\mu_k.
\]
We want to bound the second term on the right by a constant times the first plus up to the square
of the $L^2$ norm of $\varphi$, so we use the bound for $q_k$ to obtain
\begin{eqnarray*} \int_{\Sig_k}q_k\varphi^2\rho_{k+1}\ d\mu_k&\leq& c\int_{\Sig_k}(1+\sum_{p=k}^{n-1}|\tilde{A}_p|^2+\sum_{p=k+1}^{n-1}|\nabla_k\log u_p|^2)\varphi^2\rho_{k+1}d\mu_k\\
&+&\int_{\Sig_k}div_k({\cal X}_k)\varphi^2\rho_{k+1}d\mu_k.\\
\end{eqnarray*}
Now since $\varphi$ has compact support we have 
\[ \int_{\Sig_k}div_k({\cal X}_k)\varphi^2\rho_{k+1}\ d\mu_k=-\int_{\Sig_k}\langle {\cal X}_k,\nabla(\varphi^2\rho_{k+1})\rangle\ d\mu_k.
\]
Easy estimates then imply the bound
\[ |\int_{\Sig_k}div_k({\cal X}_k)\varphi^2\rho_{k+1}\ d\mu_k|\leq \frac{1}{2}
\int_{\Sig_k}|\nabla_k\varphi|^2\rho_{k+1}\ d\mu_k+c\int_{\Sig_k}(\sum_{p=k+1}^{n-1}|\nabla_k\log u_p|^2)\varphi^2\rho_{k+1}\ d\mu_k.
\]
We may now absorb the first term back to the left and use (\ref{eqn:qbound}) to obtain the bound
\[ \int_{\Sig_k}(|\nabla_k\varphi|^2+P_k\varphi^2)\rho_{k+1}\ d\mu_k\leq cQ_k(\varphi,\varphi)+
\int_{\Sig_k}\varphi^2\rho_{k+1}d\mu_k.
\]

To bound the term involving $|\nabla_klog\ u_k|^2$ we recall that on the regular set
we have
\[ \tilde{\Delta}_ku_k+q_ku_k=-\lambda_ku_k
\]
where $\lambda_k\geq 0$. This implies by direct calculation
\[ \tilde{\Delta}log\ u_k=-q_k-\lambda_k-|\nabla_klog\ u_k|^2.
\]
(Note that $\tilde{\nabla}_k=\nabla_k$ on functions which do not depend on the extra variables
$t_p$.) Now if $\varphi$ has compact support in ${\cal R}_k$, we multiply by $\varphi^2$, 
integrate by parts to obtain
\[ \int_{\Sig_k}(|\nabla_klog\ u_k|^2+q_k)\varphi^2\rho_{k+1}\ d\mu_k\leq 2\int_{\Sig_k}\varphi\langle\nabla_k\varphi, \nabla_klog\ u_k\rangle\rho_{k+1}\ d\mu.
\]
By the arithmetic-geometric mean inequality
\begin{eqnarray*}
\int_{\Sig_k}(|\nabla_klog\ u_k|^2+q_k)\varphi^2\rho_{k+1}\ d\mu_k&\leq& \frac{1}{2}\int_{\Sig_k}(|\nabla_klog\ u_k|^2+q_k)\varphi^2\rho_{k+1}\ d\mu_k \\
&+&2\int_{\Sig_k}|\nabla_k\varphi|^2\rho_{k+1}\ d\mu_k.
\end{eqnarray*}
This implies
\[ \frac{1}{2}\int_{\Sig_k}|\nabla_klog\ u_k|^2\varphi^2\rho_{k+1}\ d\mu_k\leq\frac{1}{2}Q_k(\varphi,\varphi)
+\frac{3}{2}\int_{\Sig_k}|\nabla_k\varphi|^2\rho_{k+1}\ d\mu_k.
\]
The first inequality then follows from this and our previous estimate.

The second conclusion follows since $|\nabla_k log\ \rho_{k+1}|^2\leq cP_k$, and
so the integrand on the left $|\nabla_k(\varphi\sqrt{\rho_{k+1}})|^2+|A_k|^2\varphi^2\rho_{k+1}$
is bounded pointwise by a constant times $(|\nabla_k\varphi|^2+P_k\varphi^2)\rho_{k+1}$.
\end{pf}

Recall that an important analytic step in the minimal hypersurface regularity theory is
the local reduction to the case in which the hypersurface is the boundary of a set. This
makes comparisons particularly simple and reduces consideration to a multiplicity one
setting. We will need an analogous reduction in our situation. Since the leaves of a $k$-slicing 
can be singular, we must consider the possibility that local topology comes into play and
prohibits such a reduction to boundaries of sets. What saves us here is the fact that $k$-slicings
come with a natural trivialization of the normal bundle (on the regular set). We have the following
result.
\begin{prop}
\label{prop:boundary} Assume that $U$ is compactly contained in $\O$, and that $U\cap\Sig_n$ is diffeomorphic to a ball. Assume that we have a minimal $k$-slicing in $\O$ such that the associated
$(k+1)$-slicing is partially regular. Let $\hat{\Sig}_k$ denote the closure of any connected
component of $\Sig_k\cap U\cap{\cal R}_{k+1}$. Then it follows that $\hat{\Sig}_k$ divides
the corresponding connected component (denoted $\hat{\Sig}_{k+1}$) of $\Sig_{k+1}$ into a union of two relatively open subsets, and choosing the one, denoted $U_{k+1}$,
for which the unit normal of $\hat{\Sig}_k$ points outward, we have $\hat{\Sig}_k=\partial U_{k+1}$
as a point set boundary in $\hat{\Sig}_{k+1}$, and as an oriented boundary in ${\cal R}_{k+1}$.  
\end{prop}
\begin{pf} Since $\hat{\Sig}_k\cap{\cal R}_{k+1}$ and $\hat{\Sig}_{k+1}\cap{\cal R}_{k+1}$ are
connected, it follows that the complement of $\hat{\Sig}_k\cap{\cal R}_{k+1}$ in 
$\hat{\Sig}_{k+1}\cap{\cal R}_{k+1}$ has either $1$ or $2$ connected components. These consist
of the connected components of points lying near $\hat{\Sig}_k$ on either side. Locally these are
separate components, but they may reduce globally to a single connected component. If this were
to happen, then since $dim({\cal S}_{k+1})\leq k-2$, we could find a smooth embedded closed curve $\gamma(t)$ parametrized
by a periodic variable $t\in [0,1]$ with $\gamma(0)\in\hat{\Sig}_k\cap{\cal R}_{k+1}$ and
$\gamma(t)\in {\cal R}_{k+1}\sim \hat{\Sig}_k$ for $t\neq 0$. We may also assume that
$\gamma'(0)$ is transverse to $\hat{\Sig}_k$. We choose local coordinates $x^1,\ldots, x^k$
for $\hat{\Sig}_k$ in a neighborhood $V$ of $\gamma(0)$ and we may find an embedding
$F$ of $V\times S^1$ in ${\cal R}_{k+1}$ with the property that $F(0,t)=\gamma(t)$, $F(x,0)\in \hat{\Sig}_k$, $F(x,t)\not\in \hat{\Sig}_k$ for $t\neq 0$, and $\frac{\partial F}{\partial t}(x,0)$ is transverse
to $\hat{\Sig}_k$. The $k$-form $\o=\z(x)dx^1\wedge\ldots\wedge dx^k$, where $\z$
is a nonnegative and nonzero function with compact support in $V$, is a closed form which has positive integral over
$\hat{\Sig}_k$. Since the image $V_1=F(V\times S^1)$ is compactly contained in ${\cal R}_{k+1}$ and the normal bundle of $\hat{\Sig}_{k+1}$ is trivial, we may choose coordinates $x^{k+2},\ldots, x^n$ for a normal disk, and the coordinates $x^1,\ldots, x^k,t,x^{k+2},\ldots, x^n$ are then coordinates on a
neighborhood of $V_1$ in $\Sig_n$. We may then extend $\o$ to an $(n-1)$-form on this neighborhood by setting
\[ \o_1=\o\wedge \z_1(x^{k+2},\ldots,x^n)dx^{k+2}\wedge\ldots\wedge dx^n
\]
where $\z_1$ is a nonzero, nonnegative function with compact support in the domain of
$x^{k+1},\ldots,x^n$. Thus $\o_1$ is a closed $(n-1)$-form with compact support in $U\cap \Sig_n$
which has positive integral on $\hat{\Sig}_{n-1}$, the connected component of $\Sig_{n-1}$
containing $\gamma(0)$. This contradicts the condition that each connected component
of $\Sig_{n-1}$ must divide the ball $U\cap \Sig_n$ into $2$ connected components and
is the oriented boundary of one of them, say $\hat{\Sig}_{n-1}=\partial U_n$, since
Stokes theorem would imply that $\int_{\hat{\Sig}_{n-1}}\o_1=\int_{U_n}d\o_1=0$ (note that
$\o_1$ has compact support in $U\cap\Sig_n$). 
\end{pf}

We will prove a boundedness theorem which will be needed in the proof of the
compactness theorem. Note that we will obtain the partial regularity theorem by finite
induction down from dimension $n-1$, so we may assume in the following theorems that
we have already established partial regularity for $(k+1)$-slicings. In the following result
we will consider the restriction of a $k$-slicing to a small ball $B_\s(x)$ where $x\in \R^N$.
We consider the rescaled $k$-slicing of the unit ball given by $\Sig_{j,\s}=\s^{-1}(\Sig_j-x)$
with $u_{j,\s}(y)=a_ju_j(x+\s y)$ with $a_j$ chosen so that $\int_{\Sig_{j,\s}}(u_{j,\s})^2\rho_{j+1,\s}\ d\mu_j=1$. We note that by Proposition \ref{prop:boundary} we may assume that each $\Sig_j$ in $B_\s(x)$ is the oriented boundary
of a relatively open set $O_{j+1}\subseteq \Sig_{j+1}$. We take $O_{j+1,\s}$ to be the rescaled open set.
The following result implies that the rescaled $k$-slicing remains $\Lambda$-bounded for
a suitably chosen $\Lambda$. 
\begin{thm}
\label{thm:bdness} Assume that all bounded $(k+1)$-slicings are partially regular. If 
we take any bounded minimal $k$-slicing $(\Sig_j,u_j)$ in $\Omega$ and a ball $B_\s(x)$
compactly contained in $\Omega$, then there is a $\Lambda$ depending only on $\Sig_n$ such that 
$(\Sig_{j,\s},u_{j,\s})$, $j=k,\ldots,n-1$ is $\Lambda$-bounded in $B_{1/2}(0)$.
\end{thm}
\begin{pf} The proof is by a finite induction beginning with $k=n-1$. The boundedness of
$\mu_{n-1}(\Sig_{n-1,\s})$ follows by comparison with a portion of the sphere of radius $1$ 
in a standard way (see a similar argument below). We normalize 
$\int_{\Sig_{n-1,\s}}(u_{n-1,\s})^2\ d\mu_{n-1}=1$, so it remains to show
\[ \int_{\Sig_{n-1,\s}\cap B_{1/2}(0)}|A_{n-1,\s}|^2u_{n-1,\s}^2\ d\mu_{n-1}\leq \Lambda.
\]  
To see this, we use stability with the variation $\zeta u_{n-1,\s}$ to obtain
\[ \frac{1}{4}\int_{\Sig_{n-1,\s}}|A_{n-1,\s}|^2\zeta^2u_{n-1,\s}^2\ d\mu_{n-1}\leq 
Q_{n-1,\s}(\zeta u_{n-1,\s},\zeta u_{n-1,\s}).
\]
Now we have by direct calculation for any $W_{1,2}(\Sig_{n-1,\s})$ function $v$
\[ Q_{n-1,\s}(\zeta v,\zeta v)=Q_{n-1,\s}(\zeta^2 v,v)+
\int_{\Sig_{n-1,\s}}v^2|\nabla_{n-1,\s}\zeta|^2\ d\mu_{n-1}.
\]
Taking $v=u_{n-1,\s}$ and choosing $\zeta$ to be a function which is $1$ on $B_{1/2}(0)$ with 
support in $B_1(0)$ and with bounded gradient we find
\[ \int_{\Sig_{n-1,\s}}|A_{n-1,\s}|^2u_{n-1,\s}^2\ d\mu_{n-1}\leq 4\lambda_{n-1,\s}+c\leq \Lambda
\]
for a constant $\Lambda$ where we have used the eigenvalue condition 
\[ Q_{n-1,\s}(\zeta^2 u_{n-1,\s},u_{n-1,\s})=\lambda_{n-1,\s}\int_{\Sig_{n-1,\s}}\zeta^2u_{n-1,\s}^2\ d\mu_{n-1}
\]
and the obvious relation $\lambda_{n-1,\s}=\s^2\lambda_{n-1}$. This proves $\Lambda$-boundedness for $k=n-1$.

Now assume that we have $\Lambda$-boundedness for $j\geq k+1$ in $B_{3/4}(0)$.  Thus it follows
that $\int_{\Sig_{k+1,\s}\cap B_{3/4}(0)}(1+(u_{k+1,\s})^2)\rho_{k+2,\s}\ d\mu_{k+1}$ is bounded
and hence $\int_{\Sig_{k+1,\s}\cap B_{3/4}(0)}\rho_{k+1,\s}\ d\mu_{k+1}$ is bounded. We may then
use the coarea formula to find a radius $r\in (1/2,3/4)$ so that
\[ \int_{\Sig_{k+1,\s}\cap \partial B_r(0)}\rho_{k+1,\s}\ d\mu_k\leq \Lambda.
\]
Using the portion of $\Sig_{k+1,\s}\cap \partial B_r(0)$ lying outside $O_{k,\s}$ as a comparison
surface we find 
\[ Vol_{\rho_{k+1,\s}}(\Sig_{k,\s}\cap B_{1/2}(0))\leq Vol_{\rho_{k+1,\s}}(\Sig_{k+1,\s}\cap \p B_r(0))\leq
\Lambda.
\]
Finally we prove the bound
\[ \int_{\Sig_{k,\s}\cap B_{1/2}(0)}(|A_{k,\s}|^2+\sum_{p=k+1}^n|\nabla_{k,\s} log\ u_{p,\s}|^2)u_{k,\s}^2
\rho_{k+1,\s}\ d\mu_k\leq \Lambda
\]
by the use of stability as we did above for the case $k=n-1$.
\end{pf}

We will now formulate and prove a compactness theorem for minimal $k$-slicings under
the assumption that the associated
$(k+1)$-slicings for the sequence are partially regular. We will say that a $\Lambda$-bounded sequence of $k$-slicings $(\Sig_j^{(i)},u_j^{(i)})$, $j=k,\ldots, n-1$ {\it converges} to a minimal
$k$-slicing $(\Sig_j,u_j)$ in an open set $U$ if $\Sig_j^{(i)}$ converges in $C^2$ norm to $\Sig_j$ in $\bar{U}$ locally on the complement of the singular set (of the limit)  ${\S}_j$, and such that for 
$j=k,\ldots,n-1$
\begin{equation}
\label{eqn:conv1} \lim_{i\to\infty}V_{\rho_{j+1}^{(i)}}(\Sig_j^{(i)}\cap U_i)=V_{\rho_{j+1}}(\Sig_j\cap U),
\end{equation}
\begin{eqnarray}
\label{eqn:conv2}
\lim_{i\to\infty}\|u^{(i)}_j\|^2_{0,j,U_i}&=&\|u_j\|_{0,j,U}^2 \\
\lim_{i\to\infty}  \int_{\Sig_j^{(i)}\cap U_i}(|\nabla_ju_j^{(i)}|^2+P_j^{(i)}(u_j^{(i)})^2)\rho_{j+1}^{(i)}\ d\mu_j&=&\int_{\Sig_j\cap U}(|\nabla_ju_j|^2+P_ju_j^2)\rho_{j+1}\ d\mu_j \nonumber
\end{eqnarray} 
where $U_i$ is a sequence of compact subdomains of $U$ with $U_i\subseteq U_{i+1}\subseteq U$ 
and $U=\cup_iU_i$.

To make precise the meaning of convergence on compact subsets for this
problem involves some subtlety since changing the $u_p$, $p\geq j+1$ by multiplication by a positive constant has no effect on the $\Sig_j$, so in order to get nontrivial limits for the $u_p$ we must normalize
them appropriately. In case $\Sig_j\cap U$ has multiple components this normalization must be
done on each component. If $(\Sig_j,u_j)$ is a minimal $k$-slicing with $\Sig_j$ being partially regular
for $j\geq k+1$, then we call a compact subdomain $U$ of $\O$ {\it admissible for $(\Sig_j,u_j)$} 
if $U$ is a smooth domain which meets $\partial\Sig_j$ transversally and $dim(\partial U\cap{\cal S}_j)\leq j-3$. It follows from the coarea formula that any smooth domain can be
perturbed to be admissble. We make the following definition.
\begin{defn} We say that a sequence of $k$-slicings $(\Sig_j^{(i)},u_j^{(i)})$ {\it converges 
on compact subsets} to a $k$-slicing $(\Sig_j,u_j)$ if for any compact subdomain $U$ of $\O$ which
is admissible for $(\Sig_j,u_j)$ and for any admissible domains $U_i$ for $(\Sig_j^{(i)},u_j^{(i)})$
with $U_i\subseteq U_{i+1}\subseteq U$ compactly contained in $U$ it is true that each connected component of  $\Sig_j\cap{\cal R}_{j+1}\cap U$ is a limit of 
connected components of $\Sig_j^{(i)}\cap{\cal R}_{j+1}^{(i)}\cap U_i$ in the sense of (\ref{eqn:conv1}) and (\ref{eqn:conv2}) with $u_j$ appropriately normalized on each connected component.
\end{defn}
\begin{rem} Because of the connectedness of the regular set and the Harnack inequality, we may normalize the $u_j$ to be equal to $1$ at a point of $x_0\in {\cal R}_k$ about which we have a uniform ball on which the $\Sig_j$ have bounded curvature, and this normalization suffices for the connected component of $\Sig_k\cap U$ for any compact admissible domain for $(\Sig_j,u_j)$. A consequence of
the compactness theorem below implies that this normalization suffices.
\end{rem}

The following compactness and regularity theorem includes Theorem \ref{thm:reg} as
a special case. 
\begin{thm}
\label{thm:cptness} Assume that all bounded minimal $(k+1)$-slicings are partially regular.
Given a $\Lambda$-bounded sequence of $k$-slicings , there is a subsequence which 
converges to a $\Lambda$-bounded $k$-slicing on compact open subsets of $\Omega$. 
Furthermore $\Sig_k$ is partially regular.
\end{thm}
\begin{pf} We will proceed as usual by downward induction beginning with $k=n-1$. We
will break the proof into two separate steps, the first establishing the first statement
of (\ref{eqn:conv1}) for convergence of the $\Sig_k$ and the second showing the other
two statements (\ref{eqn:conv2}) involving convergence of the $u_k$. For $k=n-1$ the first step follows 
from the usual compactness theorem for volume minimizing hypersurfaces (see \cite{simon}).
To complete the proof we will need to develop some monotonicity ideas both for the 
$\Sig_j$ and for the $u_j$. We digress on this topic and return to the proof below.

We now prove a version of the monotonicity of the frequency-type function. This idea
is due to F. Almgren \cite{almgren}, and it gives a method to prove that solutions
of variationally defined elliptic equations are approximately homogeneous on a small scale.
The importance of this method for us is that it works in the presence of singularites provided
certain integrals are defined.
We will apply this to show that the $u_k$ become homogeneous upon rescaling at a given
singular point. Assume that $C$ is a $k$ dimensional cone in $\R^n$ which is regular
except for a set $\S$ with $dim({\S})\leq k-3$. Assume that $Q$ is a quadratic form on $C$ of
the form
\[ Q(\varphi,\varphi)=\int_C(|\nabla\varphi|^2-q(x)\varphi^2)\rho\ d\mu
\]
where $\rho$ is a homogeneous weight function on $C$ of degree $p$; i.e. assume that
$\rho(\lambda x)=\lambda^p\rho(x)$ for $x\in C$ and $\lambda>0$. Assume also that
$\rho$ is smooth and positive on the regular set ${\cal R}$ of $C$ and that $\rho$ is locally 
$L^1$ on $C$. Assume also that $q$ is smooth on ${\cal R}$ and is homogeneous of 
degree $-2$; i.e. assume that
$q(\lambda x)=\lambda^{-2}q(x)$ for $x\in C$ and $\lambda>0$. Finally assume that
$u$ is a minimizer for $Q$ in a neighorhood of $0$ and in particluar that $u$ is smooth and
positive on ${\cal R}$. Assume also that $q=div({\cal X})+\bar{q}$ where $|{\cal X}|^2+|\bar{q}|\leq P$
for some positive function $P$ and that the following integral bound holds
\[ \int_C[|\nabla u|^2+(1+|\nabla log\ \rho|^2+P)u^2]\rho\ d\mu<\infty. 
\]
Under these conditions we may define
the frequency function $N(\sigma)$ which is a function of a radius $\s>0$ such that $B_\s(0)$
is contained in the domain of definition of $u$. It is defined by
\begin{equation}
\label{eqn:freqfcn} N(\s)=\frac{\s Q_{\s}(u)}{I_{\s}(u)}
\end{equation} 
where $Q_\s(u)$ and $I_\s(u)$ are defined by 
\[ Q_\s(u)=\int_{C\cap B_\s(0)}(|\nabla u|^2-q(x)u^2)\rho\ d\mu_k,\ I_\s(u)=\int_{C\cap \partial B_\s(0)} u^2\rho\ d\mu_{k-1}
\]
where the last integral is taken with respect to $k-1$ dimensional Hausdorff measure. We
may now prove the following monotonicity result for $N(\s)$.
\begin{thm}
\label{thm:freq} Assume that $u$ is a critical point of $Q$ which is integrable as above. The function $N(\s)$ is monotone increasing in $\s$, and for almost all $\s$ we have
\[ N'(\s)=\frac{2\s}{I_\s(u)}(I_\s(u_r)I_\s(u)-\langle u_r, u\rangle_\s^2)
\] 
where $u_r$ denotes the radial derivative of $u$ and $\langle\cdot,\cdot\rangle_\s$ denotes
the $\rho$-weighted $L^2$ inner product taken on $C\cap\partial B_\s(0)$. The limit of $N(\s)$
as $\s$ goes to $0$ exists and is finite. The function
$N(\s)$ is equal to a constant $N(0)$ if and only if $u$ is homogeneous of degree $N(0)$.
\end{thm}
\begin{pf} The argument can be done variationally and combines two distinct deformations of
the function $u$. The first involves a radial deformation of $C$; precisely, let $\zeta(r)$ be
a function which is nonnegative, decreasing, and has support in $B_\s(0)$. Let $X$ denote
the vector field on $\R^n$ given by $X=\z(r) x$ where $x$ denotes the position vector. The flow 
$F_t$ of $X$ then preserves $C$, and we may write
\[ Q_\s(u\circ F_t)=\int_{C\cap B_\s(0)}(|\nabla_t u|^2-(q\circ F_{t})u^2)\rho\circ F_{t}\ d\mu_t 
\]
where we have used a change of variable and $\nabla_t$ and $\mu_t$ denotes the gradient
operator and volume measure with respect to $F_t^*(g)$ where $g$ is the induced metric on $C$
from $\R^n$. Differentiating with respect to $t$ and setting $t=0$ we obtain
\[ 0=\int_C\{(\langle-{\cal L}_Xg,du\otimes du\rangle-X(q)u^2)\rho+(|\nabla u|^2-qu^2)
(X(\rho)+\rho\ div(X))\}\ d\mu
\]
where ${\cal L}$ denotes the Lie derivative. By direct calculation we have $X(q)=-2\z q$,
$X(\rho)=p\z \rho$, $div(X)=r\z'(r)+k\z$, and ${\cal L}_X g=2r\z'(r)(dr\otimes dr)+2\z g$. 
Substituting in this information and collecting terms we have
\[ 0=\int_C\{(p+k-2)\z(|\nabla u|^2-qu^2)+r\z'(|\nabla u|^2-2u_r^2-qu^2)\}\ \rho\ d\mu.
\]
Letting $\z$ approach the characteristic function of $B_\s(0)$ this implies
\begin{eqnarray*}
(p+k-2)Q_\s(u)&=&\s\int_{C\cap\partial B_\s(0)}(|\nabla u|^2-2u_r^2-qu^2)\}\ \rho\ d\mu_{k-1}\\
&=&\s\frac{dQ_\s(u)}{d\s}-2\s\int_{C\cap\partial B_\s(0)}u_r^2\rho\ d\mu_{k-1}.
\end{eqnarray*}

The second ingredient we need comes from the deformation $u_t=(1+t\z(r))u$ where
$\z$ is as above. Since $\dot{u}=\z u$ this deformation implies
\[ 0=\int_C(\langle \nabla u,\nabla(\z u)\rangle-q\z u^2)\rho\ d\mu.
\]
Expanding this and letting $\z$ approach the characteristic function of $B_\s(0)$ we have
\[ Q_\s(u)=\int_{C\cap\partial B_\s(0)}uu_r\ \rho\ d\mu_{k-1}.
\]

The proof will now follow by combining these. First we have
\[ N'(\s)=I_\s(u)^{-2}\{(Q_\s+\s Q'_\s)I_\s-\s Q_\s I'_\s\}.
\]
Substituting in for the terms involving derivatives this implies
\begin{eqnarray*} N'(\s)&=&I_\s^{-2}\{(Q_\s+(p+k-2)Q_\s)I_\s
-Q_\s(p+k-1)I_\s)\} \\
&+&2\s I_\s^{-2}\{\int_{C\cap\partial B_\s(0)} u_r^2\rho\ d\mu_{k-1}-Q_\s^2I_\s\}. 
\end{eqnarray*}
Since the first term on the right is $0$, we may write this as
\[ N'(\s)=2I_\s(u)^{-1}(I_\s(u)I_\s(u_r)-\langle u_r, u\rangle_\s^2)
\]
which is the desired formula.

To see that $N(\s)$ is bounded from below as $\s$ goes to $0$ we can observe that
\[ N(\s)=\frac{1}{2}\s\frac{d}{d\s}\log(\bar{I}_\s(u)),\ \bar{I}_\s(u)=\frac{\int_{C\cap\partial B_\s(0)}u^2\rho\ d\mu_{k-1}}{\int_{C\cap\partial B_\s(0)}\rho\ d\mu_{k-1}},
\]
and the monotonicity expresses the condition that the function $\log\ \bar{I}_\s(u)$ is a convex function
of $t=\log \s$. Since this function is defined for all $t\leq 0$, and by the coarea formula for any $\s_1>0$, there is a $\s\in [\s_1,2\s_1]$ so that $I_\s(u)\leq c\s^{-1}$ it follows that there is a sequence 
$t_i=\log\ \s_i$
tending to $-\infty$ such that $\bar{I}_{\s_i}(u)\leq c\s_i^{-K}$ for some $K>0$. Thus we have the
function $\log\ \bar{I}_{\s_i}(u)\leq -ct_i$. It follows that the slope (that is $N(\s)$) of the convex function
$\log\bar{I}_\s(u)$ is bounded from below as $t$ tends to $-\infty$.

Now if $N(\s)=N(0)$ is constant, we must have equality in the Schwartz inequality for
each $\s$, and hence we would have $u_r=f(r)u$ for some function $f(r)$. Now
this implies that $Q_\s=f(\s)I_\s$ and hence we have $rf(r)=N(0)$. Therefore it follows
that $f(r)=r^{-1}N(0)$, and $ru_r=N(0)u$ so $u$ is homogeneous of degree $N(0)$ by
Euler's formula.
\end{pf}
We will need to extend the usual monotonicity formula for the volume of minimal
submanifolds to the setting in which the submanifold under consideration minimizes
a weighted volume with a homogeneous weight function within a partially regular cone.
Precisely, let $C$ be a $k+1$ dimensional cone in $\R^n$ with a singular set $\S$ of Hausdorff
dimension at most $k-2$. Let $\rho$ be a positive weight function which is homogeneous
of degree $p$; i.e. we have $\rho(\lambda x)=\lambda^p\rho(x)$ for $x\in C$ and $\lambda>0$. 
Assume that $\rho$ is smooth and positive on the regular set of $C$, and that $\rho$ is locally
integrable with respect to Hausdorff measure on $C$. 
\begin{thm}
\label{thm:mncty} Let $\Sigma$ be a hypersurface in a $k+1$ dimensional cone $C$ which
minimizes the weighted volume $V_\rho$ for a homogeneous weight function $\rho$. We then
have the monotonicity formula
\[ \frac{d}{d\s}(\s^{-k-p}Vol_\rho(\Sig\cap B_\s(0))=
\int_{\Sig\cap\partial B_\s(0)}r^{-p-k-2}|x^\perp|^2\rho\ d\mu_{k-1}
\] 
where $x^\perp$ denotes the component of the position vector $x$ perpendicular to $\Sig$.
\end{thm}
\begin{pf} We take a function $\z(r)$ which is decreasing, nonnegative, and equal to $0$
for $r>\s$, and we consider the vector field $X=\z x$ where $x$ denotes the position vector.
The first variation formula for the $\rho$-weighted volume then implies
\[ 0=\int_\Sig(X(\rho)+div_\Sig(X)\rho)\ d\mu_k.
\]
Since $\rho$ is homogeneous we have $X(\rho)=p\z \rho$, and by direct
calculation $div_\Sig(X)=k\z+r^{-1}\z'|x^T|^2$ where $x^T$ denotes the component of $x$ tangential
to $\Sig$. Thus we have
\[ 0=\int_\Sig\{(p+k)\z+r^{-1}\z'|x^T|^2\}\rho\ d\mu_k
\]
Taking $\z$ to approximate the characteristic function of $B_\s(0)$ we may write this
\[ (p+k)Vol_\rho(\Sig\cap B_\s(0))=\s\frac{d}{d\s}Vol_\rho(\Sig\cap B_\s(0))-
\int_{\Sig\cap\partial B_\s(0)} r^{-1}|x^\perp|^2\rho\ d\mu_{k-1}
\]
where $x^\perp$ is the component of $x$ normal to $\Sig$ in $C$. Note that 
$r^2=|x^T|^2+|x^\perp|^2$ because $C$ is a cone and so $x$ is tangential to $C$. This may
be rewritten as the desired monotonicity formula and completes the proof.
\end{pf}

We now show that there can be no tangent minimal $2$-slicing with $C_2$ having an
isolated singularity at $\{0\}$.
\begin{thm}
\label{thm:2dcone} If $C_2$ is a cone lying in a tangent minimal $2$-slicing such that
$C_2\sim\{0\}\subseteq{\cal R}_2$, then $C_2$ is a plane and ${\cal R}_2=C_2$. 
\end{thm}
\begin{pf} From the eigenvalue estimate of Theorem \ref{thm:eval} we have
\[ \int_{C_2}(\frac{3}{4}\sum_{j=3}^n|\nabla_2\ log\ u_j|^2-R_2)\varphi^2\ d\mu_2\leq 4\int_{C_2}|\nabla_2\varphi|^2\ d\mu_2
\]
for test functions $\varphi$ with compact support in $C_2\sim \{0\}$. Since $C_2$ is a two
dimensional cone we have $R_2=0$ away from the origin, and hence we have
\[ \int_{C_2}\sum_{j=3}^n|\nabla_2\ log\ u_j|^2\varphi^2\ d\mu_2\leq c\int_{C_2}|\nabla_2\varphi|^2\ d\mu_2.
\]
Letting $r$ denote the distance to the origin, we take $\e$ and $R$ so that $0<\e<<R$ and
choose $\varphi$ to be a function of $r$ which is equal to $0$ for $r\leq \e^2$, equal to $1$
for $\e\leq r\leq R$, and equal to $0$ for $r\geq R^2$. In the range $\e^2\leq r\leq \e$ we 
choose 
\[ \varphi(r)=\frac{log(\e^{-2}r)}{log(\e^{-1})}
\]
and for $R\leq r\leq R^2$
\[ \varphi(r)=\frac{log(R^2r^{-1})}{log\ R}.
\] 
Thus for $\e^2\leq r\leq \e$ we have $|\nabla_2\varphi|^2=(r|log\ \e|)^{-2}$ and for
$R\leq r\leq R^2$ we have $|\nabla_2\varphi|^2=(r\ log\ R)^{-2}$. It thus follows
that
\[ \int_{C_2}|\nabla_2\varphi|^2\ d\mu_2\leq c(|log\ \e|^{-1}+(log\ R)^{-1}).
\]
Thus we may let $\e$ tend to $0$ and $R$ tend to $\infty$ to conclude that the functions
$u_3,\ldots, u_n$ are constant on $C_2$. This implies that $C_2$ has zero mean curvature and
hence is a plane. If all of the cones $C_3,\ldots C_{n-1}$ are regular near the origin, then it
follows that $0\in {\cal R}_2$, and we have completed the proof. Otherwise there is a $C_m$ for 
$m\geq 3$ which denotes the largest dimensional cone in the minimal $2$-slicing for which the 
origin is a singular point. It follows that $C_m$ is a volume minimizing cone in 
$\R^{m+1}=C_{m+1}$, and hence $u_m$ must be homogeneous of a negative degree (see Lemma \ref{lem:negdeg} below) contradicting
the fact that $u_m$ is constant along $C_2$. This completes the proof.
\end{pf}
{\it Completion of proof of Theorem \ref{thm:cptness}:} We first prove the compactness of the
$\Sig_k$ in the sense of (\ref{eqn:conv1}) under the assumption that we have the partial regularity of 
bounded minimal $(k+1)$-slicings and the compactness (both (\ref{eqn:conv1}) and (\ref{eqn:conv2}))
for $j\geq k+1$. We need the following lemma.
\begin{lem} 
\label{lem:locbd} Assume that both the compactness and partial regularity hold for $(k+1)$-slicings. Given any $x\in {\S}_{k+1}$, there are constants $c$ and $r_0$ (depending on $x$ and $\Sig_{k+1}$) 
so that for $r\in (0,r_0]$ we have
\[  \int_{\Sig_{k+1}\cap B_{2r}(x)} u_{k+1}^2\rho_{k+2}\ d\mu_{k+1}\leq 
cr^2\int_{\Sig_{k+1}\cap B_r(x)}P_{k+1}u_{k+1}^2\rho_{k+2}\ d\mu_{k+1},
\]
and
\[ Vol_{\rho_{k+2}}(\Sig_{k+1}\cap B_{2r}(x))\leq cVol_{\rho_{k+2}}(\Sig_{k+1}\cap B_r(x)).
\]
\end{lem}
\begin{pf} Since the left hand side of the inequality is continuous under convergence and the
right hand side is lower semicontinuous (Fatou's theorem) it is enough to establish the inequality
for $r=1$ on a cone $C_{k+1}$. This we can do by a compactness argument since we can
normalize 
\[ \int_{C_{k+1}\cap B_1(0)} u_{k+1}^2\rho_{k+2}\ d\mu_{k+1}=1
\] 
and if we had a sequence
of singular cones for which the right hand side tends to zero we would have a limiting cone $C_{k+1}$
on which $P_{k+1}=0$. It follows that $u_{k+2},\ldots, u_{n-1}$ are constant on $C_{k+1}$. Note
that the highest dimensional {\it singular} cone in the slicing $C_{n_0}$ is minimal and hence $u_{n_0}$ is homogeneous of a negative degree (see Lemma \ref{lem:negdeg} below). Therefore if 
$n_0>k+1$ we have a contradiction. Therefore 
we conclude that $C_{k+1}$ is minimal and $C_{k+2},\ldots, C_{n-1}$ are planes. Thus it follows
that $\tilde{A}_{k+1}=A_{k+1}=0$ and hence $C_{k+1}$ is also a plane. Thus the cones are regular
sufficiently far out in the sequence; a contradiction. The second inequality follows easily by reduction
to cones. This proves the bounds. 
\end{pf}
Given a sequence $(\Sig_j^{(i)},u_j^{(i)})$ of $\Lambda$-bounded minimal $k$- slicings, we 
may apply the inductive assumption to obtain a subsequence (with the same notation) for 
which the corresponding sequence of $(k+1)$-slicings converges in the sense of (\ref{eqn:conv1}) 
and (\ref{eqn:conv2}). By standard compactness theorems we may assume that $\Sig_k^{(i)}$
converges on compact subsets of $\Omega\sim {\S}_{k+1}$ to a limiting submanifold $\Sig_k$
which minimizes $Vol_{\rho_k}$ (and is therefore regular outside a closed set of dimension at most
$k-7$). To establish (\ref{eqn:conv1}) we choose a neighborhood $U$ of ${\S}_{k+1}$ such
that 
\[ Vol_{\rho_{k+2}}(\Sig_{k+1}\cap \bar{U})<\e.
\]
We apply Lemma \ref{lem:locbd} and compactness to find a finite collection of points 
$x_\a\in {\S}_{k+1}$ and balls  $B_{r_\a}(x_\a)\subset U$ so that
\[  \int_{\Sig_{k+1}\cap B_{2r_\a}(x_\a)} u_{k+1}^2\rho_{k+2}\ d\mu_{k+1}< 
cr_\a^2\int_{\Sig_{k+1}\cap B_{r_\a}(x_\a)}P_{k+1}u_{k+1}^2\rho_{k+2}\ d\mu_{k+1}
\]
and
\[ Vol_{\rho_{k+2}}(\Sig_{k+1}\cap B_{2r_\a}(x_\a))< cVol_{\rho_{k+2}}(\Sig_{k+1}\cap B_{r_\a}(x_\a)).
\]
Now apply the Besicovitch covering
lemma to extract a finite number of disjoint collections ${\cal B}_\a$, $\a=1,\ldots, K$ of such balls 
whose union covers ${\S}_{k+1}$. If $V$ denotes the union of these balls, then $V$ is a neighborhood 
of ${\S}_{k+1}$, and hence for
$i$ sufficiently large we have ${\S}_{k+1}^{(i)}\subset V$. Because of convergence of the left sides
and lower semicontinuity of the right side, we have for $i$ sufficiently large
\[  \int_{\Sig_{k+1}^{(i)}\cap B_{2r_\a}(x_\a)} (u_{k+1}^{(i)})^2\rho_{k+2}^{(i)}\ d\mu_{k+1}< 
cr_\a^2\int_{\Sig_{k+1}^{(i)}\cap B_{r_\a}(x_\a)} P_{k+1}^{(i)}(u_{k+1}^{(i)})^2\rho_{k+2}^{(i)}\ d\mu_{k+1}
\]
and
\[ Vol_{\rho_{k+2}^{(i)}}(\Sig_{k+1}^{(i)}\cap B_{2r_\a}(x_\a))< cVol_{\rho_{k+2}^{(i)}}(\Sig_{k+1}^{(i)}\cap B_{r_\a}(x_\a)).
\]
By the coarea formula, for each such 
ball $B_{r_0}(x)$ we may find $s\in [r_0,2r_0]$ ($s$ depending on $i$) so that 
\[ Vol_{\rho_{k+1}^{(i)}}(\Sig_{k+1}^{(i)}\cap\partial B_s(x))\leq 2r_0^{-1}\int_{\Sig_{k+1}^{(i)}\cap B_{2r_0}}
u_{k+1}^{(i)}\rho_{k+2}^{(i)}\ d\mu_{k+1}.
\]
Using the minimizing property of $\Sig_k^{(i)}$ and simple inequalities we find
\begin{eqnarray*}
 Vol_{\rho_{k+1}^{(i)}}(\Sig_k^{(i)}\cap B_{r_0})&\leq& \e_1^{-1}\int_{\Sig_{k+1}\cap B_{2r_0}(x)}\rho_{k+2}^{(i)}\ d\mu_{k+1} \\ 
 &+&\e_1r_0^{-2}\int_{\Sig_{k+1}\cap B_{2r_0}}(u_{k+1}^{(i)})^2\rho_{k+2}^{(i)}\ d\mu_{k+1}
\end{eqnarray*}
for any $\e_1>0$. Applying the inequalities above and summing over the balls (using disjointness and
a bound on $K$) we find
\[ Vol_{\rho_{k+1}^{(i)}}(\Sig_k^{(i)}\cap V)\leq c\e_1^{-1}Vol_{\rho_{k+2}^{(i)}}(\Sig_{k+1}^{(i)}\cap \bar{U})
+c\e_1\int_{\Sig_{k+1}^{(i)}}P_{k+1}^{(i)}(u_{k+1}^{(i)})^2\rho_{k+2}^{(i)}\ d\mu_{k+1}.
\] 
For $i$ sufficiently large this implies
\[ Vol_{\rho_{k+1}^{(i)}}(\Sig_k^{(i)}\cap V)\leq c\e_1^{-1}\e+c\e_1,
\]
so that we may fix $\e_1$ sufficiently small and then choose $\e$ as small as we wish to
make the right hand side smaller than any preassigned amount. Since we have
\[ \lim_{i\to\infty}Vol_{\rho_{k+1}^{(i)}}(\Sig_k^{(i)}\sim V)=Vol_{\rho_{k+1}}(\Sig_k\sim V),
\]
we can conclude that $\lim_{i\to\infty}Vol_{\rho_{k+1}^{(i)}}(\Sig_k^{(i)})=Vol_{\rho_{k+1}}(\Sig_k)$
establishing (\ref{eqn:conv1}).

Now assume that we have established the
partial regularity of all bounded minimal $(k+1)$-slicings and that we have proven the compactness
for the $\Sig_k$ in the sense of (\ref{eqn:conv1}). We can then use the results we have obtained above together
with dimension reduction to prove partial regularity for $\Sig_k$. Precisely, we have $dim({\S}_k)\leq
k-2$, and if $dim({\S}_k)>k-3$, then we can choose a number $d$ with
\[ k-3<d<dim({\cal S}_k),
\]
and go to a point $x\in {\cal S}_k$ of density for the measure ${\cal H}^d_\infty$ (since 
${\cal H}^d_\infty({\cal S}_k)>0$). Taking successive tangent cones in the standard way and using the upper-semicontinuity
of ${\cal H}^d_\infty({\S}_k)$ we would eventually produce a minimal $2$-slicing by cones such that
$C_2\times \R^{k-2}$ has singular set with Hausdorff dimension at most $k-2$ (by partial regularity
of $(k+1)$-slicings) and greater than $k-3$. Therefore the cone $C_2$ must have an isolated singularity at the origin. This in turn contradicts Theorem \ref{thm:2dcone}. Therefore it follows that $dim({\S}_k)\leq k-3$ and $\Sig_k$ is partially regular.

The final step of the proof is to show that the compactness statement holds for the $u_k$ under
the assumption that it holds for $(\Sig_j,u_j)$ for $j\geq k+1$ and also for $\Sig_k$ (as  
established above). Assume that we have a sequence of minimal $k$-slicings such that the associated 
$(k+1)$-slicings and $\Sig_k^{(i)}$ converge on compact subsets in the sense of (\ref{eqn:conv1}) and (\ref{eqn:conv2}). We choose a compact domain $U$ which is admissble for $(\Sig_j,U_j)$ and
a nested sequence of domains $U_i$ admsisible for $(\Sig_j^{(i)},u_j^{(i)})$. We work with
a connected component of $\Sig_k\cap U$ which by abuse of notation we call by the same name $\Sig_k$. 

We may assume that the $u_k^{(i)}$ converge uniformly to $u_k$ on compact subsets of 
$\Omega\sim {\cal S}_k$ (where we can write $\Sig_k^{(i)}$ locally as a normal graph over 
$\Sig_k$ and compare corresponding values of $u_k^{(i)}$ to $u_k$). In particular, if
$W$ is a compact subdomain of $\O\cap{\cal R}_k$ we have convergence of weighted
$L^2$ norms of $u_k^{(i)}$ to the corresponding $L^2$ norm of $u_k$ on $W$. If $U$
is any compact subdomain of $\O$ and
$\eta>0$, then by Proposition \ref{prop:l2con} applied with ${\cal S}={\cal S}_k$ we can
find an open neighborhood $V$ of ${\cal S}\cap\bar{U}$ so that for $i$ sufficiently large
${\cal S}_k^{(i)}\cap\bar{U}\subset V$, and 
\[ \int_{\Sig_k^{(i)}\cap V}(u_k^{(i)})^2 \rho_{k+1}^{(i)}\ d\mu_k\leq \eta\int_{\Sig_k^{(i)}\cap\O}
[|\nabla_ku_k^{(i)}|^2+(1+P_k^{(i)})(u_k^{(i)})^2]\rho_{k+1}^{(i)}\ d\mu_k.
\]
The same inequality holds for the limit, and by the boundedness of the sequence the integral
on the right is uniformly bounded. Thus by choosing $\eta$ small enough we can make the
right hand side less than any prescribed $\e>0$. On the other hand if we take 
$W=U\setminus\bar{V}$ we then have convergence of the weighted $L^2$ norms on $W$,
so we can make the difference as small as we wish on $W$. It follows that the difference of
$L^2$ norms can be made arbitrarily small on $U$.
This completes the proof that the weighted $L^2$ integrals converge.

Completing the proof will require the construction of a proper locally Lipschitz function 
$\Psi_k$ on ${\cal R}_k$
such that $u_k|\nabla_k\Psi_k|$ is bounded in $L^2(\Sig_k)$. We give the construction of such a function in 
Proposition \ref{prop:proper} below. It also follows that we may construct a subsequence
so that $\Psi_k^{(i)}$ are uniformly close to $\Psi_k$ on compact subsets of $\R^N\sim{\cal S}_k$
for $i$ large. We can now prove the second part of the convergence (\ref{eqn:conv2}). Assume
that $U\subset U_1\subset \Omega$ are compact domains. We let $\e>0$ we may choose a neighborhood $V$ of ${\cal S}_k$ so small that $\int_{V\cap\bar{U_1}}u_k^2\rho_{k+1}\ d\mu_k<\e$. Because $\Psi_k$ is proper on ${\cal R}_k$, we may choose $\Lambda$ sufficiently large that  
$E_k(\Lambda)\subset V$ where $E_k(\Lambda)$ is the
subset of $\Sig_k$ on which $\Psi_k>\Lambda$. We now let $\gamma(t)$ be a nondecreasing 
Lipschitz function such that $\gamma(t)=0$ for $t<\Lambda$, $\gamma(t)=1$ for $t>\Lambda$,
and $\gamma'(t)\leq \Lambda^{-1}$. We let $\varphi$ be a spatial cutoff function which is $1$
on $U$, $0$ outside $U_1$, and has bounded gradient. We then have the inequality by Proposition
\ref{prop:coercive}
\[ \int_{\Sig_k^{(i)}}(|\nabla_k\psi_k^{(i)}|^2+
P_k^{(i)}(\psi_k^{(i)})^2)\rho_k^{(i)}\ d\mu_j\leq cQ_k(\psi_k^{(i)},\psi_k^{(i)})
\]
where $\psi_k^{(i)}=\varphi(\gamma\circ\Psi_k^{(i)})u_k^{(i)}$. Since the support of $\psi_k^{(i)}$ is
contained in $V$ for $i$ sufficiently large we then have
\[ \int_{\Sig_k^{(i)}}(|\nabla_k\psi_k^{(i)}|^2+P_k^{(i)}(\psi_k^{(i)})^2)\rho_{k+1}^{(i)}\ d\mu_j\leq
c\int_{\Sig_k^{(i)}\cap V}(1+\Lambda^{-2}|\nabla_k\Psi_k^{(i)}|^2)(u_k^{(i)})^2\rho_{k+1}^{(i)}\ d\mu_k.
\]
Since we have convergence of the $L^2$ norms of $u_k^{(i)}$ and boundedness of the $L^2$
norms of $u_k^{(i)}|\nabla_k\Psi_k^{(i)}|$, we then conclude that
\[ \int_{\Sig_k^{(i)}}(|\nabla_k\psi_k^{(i)}|^2+
P_k^{(i)}(\psi_k^{(i)})^2)\rho_{k+1}^{(i)}\ d\mu_j\leq c\e+c\Lambda^{-2}.
\]
If we let $V_1$ be a neighborhood of ${\cal S}_k$ such that 
$\Sig_k\cap V_1\subset E_k(3\Lambda)$, then for
$i$ sufficiently large we will have $\Sig_k^{(i)}\cap V_1\subset E_k^{(i)}(2\Lambda)$
and hence 
\[ \int_{\Sig_k^{(i)}\cap V_1}(|\nabla_ku_k^{(i)}|^2+
P_k^{(i)}(u_k^{(i)})^2)\rho_{k+1}^{(i)}\ d\mu_j\leq c\e+c\Lambda^{-2}.
\]
Since this can be made arbitrarily small, we have shown (\ref{eqn:conv2}) and completed the 
proof of Theorem \ref{thm:cptness}.
\end{pf}
We will need the following lemma concerning minimal cones $C_m\subset \R^{m+1}$.
\begin{lem}
\label{lem:negdeg} Assume that $C_m$ is a volume minimizing cone in $\R^{m+1}$ and that
$u_m$ is a positive minimizer for $Q_j$ which is homogeneous of degree $d$ on $C$. There
is a positive constant $c$ depending only on $m$ so that $d\leq -c$.
\end{lem}
\begin{pf} We write $u_m=r^dv(\xi)$ where $\xi\in S^m$. If we let $\Sig=C\cap S^m$,
then $v$ satisfies the eigenvalue equation $\Delta v+1/8|A_m|^2v=-\mu v$ where we must have
$d(d+m-2)=\mu$. This implies that $d=1/2(2-m+\sqrt{(m-2)^2+4\mu})$ or 
$d=1/2(2-m-\sqrt{(m-2)^2+4\mu})$. Since $v$ and $|\nabla v|$ are in $L^2(\Sig)$ we must 
have $\mu<0$ and this implies that $d<0$.
To prove the negative upper bound on $d$ recall that the set of volume minimizing cones
is a compact set, and we have proven the compactness theorem above for the $L^2$ norms,
so if we had a sequence $(C_m^{(i)},u_m^{(i)})$ such that $d^{(i)}$ tends to $0$ we could extract
a convergent subsequence of the $(\Sig^{(i)},v^{(i)})$ which converges to $(\Sig,v)$ where
we could normalize $\int_{\Sig^{(i)}}(v^{(i)})^2\ d\mu_{m-1}=1$ (hence $\int_\Sig v^2\ d\mu_{m-1}=1$).
Since we have smooth convergence on compact subsets of the complement of the singular set
of $\Sig$ we would then have $\Delta v+5/8|A_m|^2v=0$ and therefore we would have $\mu=0$ for 
the limiting cone, a contradiction.
\end{pf}
As the final topic of this section we construct the proper functions which were used in the 
proof of Theorem \ref{thm:cptness}. This result will also be used in the next section.
\begin{prop}
\label{prop:proper} Suppose we have a $\Lambda$-bounded minimal $k$-slicing in $\Omega$.
There exists a positive function $\Psi_k$ which is locally Lipschitz on ${\cal R}_k$ and such that
for any domain $U$ compactly contained in $\Omega$, the function $\Psi_k$ is proper on
${\cal R}_k\cap\bar{U}$. Moreover, the function $u_k|\nabla_k\Psi_k|$ is bounded in 
$L^2(\Sig_k\cap U)$ for any domain $U$ compactly contained in $\Omega$.
\end{prop}
\begin{pf} We define $\Psi_k=\max\{1,log\ u_k,log\ u_{k+1},\ldots,log\ u_{n-1}\}$ and we show that
it has the properties claimed. First note that $\Psi_k$ is locally Lipschitz on ${\cal R}_k$ since
it is the maximum of a finite number of smooth functions on ${\cal R}_k$. The bound
\[ \int_{\Sig_k\cap U}(u_k|\nabla_k\Psi_k|)^2\rho_{k+1}\ d\mu_k\leq 
\max_{k\leq j\leq n-1} \int_{\Sig_k\cap U}(u_k|\nabla_klog\ u_j|)^2\rho_{k+1}\ d\mu_k
\]
together with Proposition \ref{prop:coercive} implies the $L^2(\Sig_k)$ bound claimed on 
$\Psi_k$. (Note that we may replace $\varphi$ by $\varphi u_k$ in the first inequality of Proposition 
\ref{prop:coercive} where $\varphi$ is a cutoff function which is equal to $1$ on $U$.)

It remains to prove that $\Psi_k$ is proper on ${\cal R}_k\cap\bar{U}$. Since $\bar{U}$ 
is compact it suffices to show that for any $x_0\in {\cal S}_k\cap\bar{U}$ we have
\[ \lim_{x\to x_0}\Psi_k(x)=\infty.
\]
If we let $m\geq k$ be the largest integer such that $\Sig_m$ is singular at $x_0$, then
there is an open neighborhood $V$ of $x_0$ in which $\Sig_m$ is a volume minimizing
hypersurface in a smooth Riemannian manifold. We will show that $u_m$ tends to infinity
at $x_0$ by first showing that this is true for any homogeneous approximation of $u_m$
at $x_0$. In order to construct homogeneous approximations we need to have the
compactness theorem for this top dimensional case, but our proof of compactness used 
the result we are trying to prove, so we must find another argument for establishing
(\ref{eqn:conv2}) since (\ref{eqn:conv1}) is a standard result for volume minimizing hypersurfaces
in smooth manifolds. Our proof of the first part of (\ref{eqn:conv2}) did not require the
function $\Psi_k$, so we need only deal with the second part. First recall that $dim({\cal S}_m)\leq m-7$, so it follows from a standard result that given any $\e,\d>0$ and $a\in (0,7)$ we can find a Lipschitz function $\psi$ so that $\psi=1$ in a neighborhood of ${\cal S}_m$, $\psi(x)=0$ for points $x$
with $dist(x,{\cal S}_m)\geq \d$, and 
\[ \int_{\Sig_m\cap V}|\nabla_m\psi|^a\ d\mu_m<\e^a.
\] 
We show that
\[ \int_{\Sig_m\cap V}|\nabla_m\psi|^2u_m^2\ d\mu_m\leq c\e^2.
\] 
If we can establish this inequality,
then we can complete the proof of compactness for $k=m$ in the set $V$ as in the proof of
Theorem \ref{thm:cptness}. To establish
the inequality, we observe that the equation satisfied by $u_m$ is of the form
\[ \Delta_m u_m+5/8|A_m|^2u_m+qu_m=0
\] 
where $q$ is a bounded function (since $\Sig_m$ is volume minimizing in a smooth manifold). On the
other hand the stability implies that
\[ \int_{\Sig_m} |A_m|^2\varphi^2\ d\mu_m\leq \int_{\Sig_m} (|\nabla\varphi|^2+c\varphi^2)\ d\mu_m.
\]
We may then replace $\varphi$ by $u_m^{8/5}\varphi$ and use the equation for $u_m$ to
obtain
\[ \int_{\Sig_m}|\nabla_m(u_m)^{8/5}|^2\varphi^2\ d\mu_m\leq c\int_{\Sig_m}u_m^{16/5}(|\nabla_m\varphi|^2+\varphi^2)\ d\mu_m.
\]
We may then apply the Sobolev inequality for minimal submanifolds to conclude that $u_m$ 
satisfies
\[ \int_{\Sig_m\cap V}u_m^{\frac{16m}{5(m-2)}}\ d\mu_m\leq c.
\] 
We then apply the H\"older inequality to obtain
\[ \int_{\Sig_m\cap V}|\nabla_m\psi|^2u_m^2\ d\mu_m\leq \|\nabla_m\psi\|_{\frac{16m}{3m+10}}^2
\|u_m\|_{\frac{16m}{5(m-2)}}^2.
\]
Setting $a=\frac{16m}{3m+10}<7$ we have from above
\[ \int_{\Sig_m\cap V}|\nabla_m\psi|^2u_m^2\ d\mu_m\leq c\e^2
\]
as desired. 

Thus we have the compactness theorem for $(\Sig_m,u_m)$ in $V$ and we can construct 
tangent cones to $\Sig_m$ at $x_0$ and homogeneous approximations to $u_m$ at $x_0$. By 
Lemma \ref{lem:negdeg} any such homogeneous approximation $v_m$ has strictly negative degree $d\leq -c$ on its cone $C_m$ of definition. If we let ${\cal R}_m(C)$ denote the regular set
of $C$, then it follows that for any $\mu>1$, we have
\[ \inf_{{\cal R}_m(C)\cap B_{\alpha\sigma}(0)}v_m\geq \mu\inf_{{\cal R}_m(C)\cap B_\sigma(0)}v_m
\]
for a fixed constant $\alpha\in (0,1)$ depending on $\mu$, but independent of which cone and which homogeneous approximation we choose. Note that $\Delta_mu_m\leq cu_m$ and 
$\Delta_mv_m\leq 0$, so by the mean value inequality on volume minimizing hypersurfaces 
(see \cite{bg}) we have
\[ u_m(x)\geq cr^{-m}\int_{\Sig_m\cap B_r(x)}u_m\ d\mu_m,\ v_m(x)\geq 
cr^{-m}\int_{C_m\cap B_r(x)}v_m\ d\mu_m
\]
for any $r$ so that $B_r(x_0)$ is compactly contained in $V$. It follows that the essential infima of
both $u_m$ and $v_m$ are positive on any compact subset. We now show that there exists
$\alpha\in (0,1)$ such that 
\[ \inf_{{\cal R}_m\cap B_{\alpha\sigma}(x_0)}u_m\geq 2\inf_{{\cal R}_m\cap B_\sigma(x_0)}u_m
\]
for $\sigma$ sufficiently small. If we establish this, we have finished the proof that
$u_m$ tends to infinity at $x_0$ and hence we will have the desired properness 
conclusion for $\Psi_k$. To establish this inequality we observe that if $(\Sig_m^{(i)},u_m^{(i)})$
is a sequence converging to $(\Sig_m,u_m)$ in the sense of (\ref{eqn:conv1}) and (\ref{eqn:conv2})
and $K$ is a compact set such that ${\cal R}_m\cap K\neq \phi$ we have
\[ \inf_{{\cal R}_m\cap K}u_m\leq \liminf_{i\to\infty}\inf_{{\cal R}_m^{(i)}\cap K}u_m^{(i)}\leq 
\limsup_{i\to\infty}\inf_{{\cal R}_m^{(i)}\cap K}u_m^{(i)}\leq c\inf_{{\cal R}_m\cap K}u_m
\] 
for a fixed constant $c$. The first and second inequalities are obvious, and to get the third we
observe that for a small radius $r$ and any $x\in {\cal R}_m\cap K$ we have from above
\[ u_m(x)\geq cr^{-m}\int_{\Sig_m\cap B_r(x)}u_m\ d\mu_m,
\]
and hence for $i$ sufficiently large
\[ u_m(x)\geq cr^{-m}\int_{\Sig_m^{(i)}\cap B_r(x)}u_m^{(i)}\ d\mu_m\geq \e_0\inf_{\Sig_m^{(i)}\cap B_r(x)} u_m^{(i)}
\]
for a positive constant $\e_0$. This establishes the third inequality. The proof can now be completed
by using rescalings at $x_0$ which converge to $(C_m,v_m)$ for some cone and homogeneous
function together with the corresponding result for the homogeneous case. 
\end{pf}

\bigskip
\section{\bf Existence of minimal $k$-slicings}
The main purpose of this section is to prove Theorem \ref{thm:exst}. We begin with the construction
of the eigenfunction $u_k$ assuming the $\Sig_k$ has already been constructed and is partially
regular in the sense that $dim({\cal S}_k)\leq k-3$. We define the Hilbert spaces ${\cal H}_k$ 
and ${\cal H}_{k,0}$ as in the last section, namely, ${\cal H}_k$ (respectively ${\cal H}_{k,0}$) is the
completion in $\|\cdot\|_{0,1}$ of the Lipschitz functions with compact support in ${\cal R}_k\cap\bar{\Omega}$ (respectively  ${\cal R}_k\cap\Omega$).
In order to handle boundary effects we also assume that there is a larger domain $\O_1$ which contains $\bar{\O}$ as a compact subset and that the $k$-slicing is defined and boundaryless in $\O_1$. Note that this is automatic if $\partial\Sig_j=\phi$. Thus ${\cal H}_{k,0}$ consists of those functions in ${\cal H}_k$ with $0$ boundary data on $\Sig_k\cap\O$. The quadratic form $Q_k$ is nonnegative 
definite on the Lipschitz functions
with compact support in ${\cal R}_k\cap\Omega$, and so the standard Schwartz inequality holds
for any pair of such functions $\varphi,\psi$
\begin{equation}
\label{eqn:schwartz} Q_k(\varphi,\psi)\leq \sqrt{Q_k(\varphi,\varphi)}\sqrt{Q_k(\psi,\psi)}.
\end{equation}
We now have the following result.
\begin{thm}
\label{thm:qcomplete} The function $Q_k(\varphi,\psi)$ is continuous with respect to the norm
$\|\cdot\|_{0,1}$ in both variables and therefore extends as a continuous nonnegative definite
bilinear form on ${\cal H}_{k,0}$. The Schwartz inequality (\ref{eqn:schwartz}) holds for 
$\varphi,\psi\in {\cal H}_{k,0}$. The function $Q_k(\varphi,\varphi)$ is strongly continuous and 
weakly lower semicontinuous on ${\cal H}_{k,0}$.
\end{thm}
\begin{pf} From Proposition \ref{prop:coercive} we have for $\varphi_1,\varphi_2$ Lipschitz functions with compact support in ${\cal R}_k\cap\Omega$
\[ Q_k(\varphi_1-\varphi_2,\varphi_1-\varphi_2)\leq c\|\varphi_1-\varphi_2\|_{1,k}^2,
\]
so it follows from (\ref{eqn:schwartz}) that 
\[ |Q_k(\varphi_1,\psi)-Q_k(\varphi_2,\psi)|\leq \sqrt{Q_k(\varphi_1-\varphi_2,\varphi_1-\varphi_2)}\sqrt{Q_k(\psi,\psi)}.
\]
Combining these we see that $Q_k$ is continuous in the first slot, and since it is symmetric in
both slots. Therefore $Q_k$ extends as a continuous nonnegative definite bilinear form on
${\cal H}_{k,0}$ and the Schwartz inequality holds on ${\cal H}_{k,0}$ by continuity.

To complete the proof we must prove that $Q_k(\varphi,\varphi)$ is weakly lower semicontinuous
on ${\cal H}_{k,0}$. Note that the square norm $\|\varphi\|_{0,k}^2+Q_k(\varphi,\varphi)$ is equivalent to $\|\varphi\|_{1,k}^2$ by Proposition \ref{prop:coercive}. Therefore these have the same bounded 
linear functionals and hence determine the same weak topology on ${\cal H}_{k,0}$. Assume we have a sequence $\varphi\in {\cal H}_{k,0}$ which converges weakly to $\varphi\in{\cal H}_{k,0}$. We then have for any $\psi\in{\cal H}_{k,0}$
\[ Q_k(\varphi,\psi)=\lim_{i\to\infty}Q_k(\varphi_i,\psi).
\] 
This implies that for $i$ sufficiently large
\[ Q_k(\varphi,\varphi)=Q_k(\varphi-\varphi_i,\varphi)+Q_k(\varphi_i,\varphi)\leq \e+
\sqrt{Q_k(\varphi_i,\varphi_i)}\sqrt{Q_k(\varphi,\varphi)}
\]
for any chosen $\e>0$. It follows that
\[ Q_k(\varphi,\varphi)\leq \sqrt{Q_k(\varphi,\varphi)}\liminf_{i\to\infty}\sqrt{Q_k(\varphi_i,\varphi_i)}
\]
which implies the desired weak lower semicontinuity.
\end{pf}
In order to construct a lowest eigenfunction $u_k$ we will need the following Rellich-type
compactness theorem.
\begin{thm}
\label{thm:rellich} The inclusion of ${\cal H}_{k,0}$ into $L^2(\Sig_k)$ is compact in the
sense that any bounded sequence in ${\cal H}_{k,0}$ has a convergent subsequence
in $L^2(\Sig_k)$.
\end{thm}
\begin{pf} This statement follows from Proposition \ref{prop:l2con} and the standard
Rellich theorem.
Assume that we have a bounded sequence $\varphi_i\in {\cal H}_{k,0}$; that is,
 $\|\varphi_i\|_{1,k}^2\leq c$.
We may extend the $\varphi_i$ to $\O_1$ be taking $\varphi_i=0$ in $\O_1\sim\O$, and by the 
standard Rellich compactness theorem we may assume by extracting a subsequence
that the $\varphi_i$ converge in $L^2$ norm on compact subsets of $\bar{\O}\sim {\cal S}_k$ and
weakly in ${\cal H}_{k,0}$ to a limit $\varphi\in {\cal H}_{k,0}$. We show that $\varphi_i$ converges
to $\varphi$ in $L^2(\Sig_k)$.
Given any $\e_1>0$, we can choose $\e>0,\ \d>0$ in Proposition \ref{prop:l2con} so that for each
$i$ we have
\[ (\int_{\Sig_k\cap V}\varphi_i^2\rho_{k+1}\ d\mu_k)^{1/2}\leq \e_1/3
\]
where $V$ is an open neighborhood of ${\cal S}_k\cap\bar{\O}$. The Fatou theorem then implies
\[ (\int_{\Sig_k\cap V}\varphi^2\rho_{k+1}\ d\mu_k)^{1/2}\leq \e_1/3
\]
Since $K=(\Sig_k\sim V)\cap\bar{\O}$
is a compact subset of $\bar{\O}\sim {\cal S}_k$, we have for $i$ sufficiently large
\[ (\int_K(\varphi_i-\varphi)^2\rho_{k+1}\ d\mu_k)^{1/2}\leq \e_1/3.
\]
Combining these bounds we find
\[ \|\varphi_i-\varphi\|_0\leq (\int_K(\varphi_i-\varphi)^2\rho_{k+1}\ d\mu_k)^{1/2}+
(\int_{\Sig_k\cap V}(\varphi_i-\varphi)^2\rho_{k+1}\ d\mu_k)^{1/2}\leq \e_1
\]
for $i$ sufficiently large. This completes the proof.
\end{pf}
We are now ready to prove the existence, positivity, and uniqueness of $u_k$ on $\Sig_k\cap\O$.
\begin{thm}
\label{thm:spectrum} The quadratic form $Q_k$ on ${\cal H}_{k,0}$ has discrete spectrum with 
respect to the $L^2(\Sig_k)$
inner product and may be diagonalized in an orthonormal basis for $L^2(\Sig_k)$. The eigenfunctions
are smooth on ${\cal R}_k\cap\O$, and if we choose a first eigenfunction $u_k$, then $u_k$
is nonzero on ${\cal R}_k\cap\O$ and is therefore either strictly positive or strictly negative
since ${\cal R}_k\cap\O$ is connected. Furthermore any first eigenfunction is a multiple of $u_k$ which
we may take to be positive.
\end{thm}
\begin{pf} This follows from the standard minmax variational procedure for defining eigenvalues
and constructing eigenfunctions. For example, to construct the lowest eigenvalue and eigenfunction
we let
\[ \lambda_k=\inf\{Q_k(\varphi,\varphi):\ \varphi\in {\cal H}_{k,0},\ \|\varphi\|_{0,k}=1\}.
\]
By Theorem \ref{thm:rellich} and Theorem \ref{thm:qcomplete} we may achieve this infimum
with a function $u_k\in{\cal H}_{k,0}$ with $\|u_k\|_{0,k}=1$. The Euler-Lagrange equation
for $u_k$ is then the eigenfunction equation with eigenvalue $\lambda_k$. The higher eigenvalues
and eigenfunctions can be constructed by imposing orthogonality constraints with respect
the $L^2(\Sig_k)$ inner product. We omit the standard details. The smoothness on ${\cal R}_k\cap\O$
follows from elliptic regularity theory.

The fact that a lowest eigenfunction $u$ is nonzero follows from
the fact that if $u\in {\cal H}_{k,0}$ then $|u|\in {\cal H}_{k,0}$ and $Q_k(u,u)=Q_k(|u|,|u|)$ a property
which can be easily checked on the dense subspace of Lipschitz functions with compact support
in ${\cal R}_k\cap\O$ and then follows by continuity. The multiplicity one property of the lowest
eigenspace follows from this property in the usual way. We omit the details.
\end{pf}

We now come to the existence results. We first discuss Theorem \ref{thm:exst}
and we then generalize the existence proof to a more precise form. Suppose $X$ is
a closed $k$-dimensional oriented manifold with $k<n$. We assume that $\Sig_n$
is a closed oriented $n$-manifold and that there is a smooth map $F:\Sig_n\to X\times T^{n-k}$
of degree $s\neq 0$. We let $\O$ denote a (unit volume) volume form of $X$ and let $\Theta=F^*\O$ 
so that $\Theta$ is a closed $k$-form on $\Sig_n$. We let $t^p$ for $p=k+1,\ldots, n$
denote the coordinates on the circles and we assume they are periodic with period $1$.
For $p=k+1,\ldots, n$ we let $\o^p$ be the closed $1$-form $\o^p=F^*(dt^p)$. The
assumption on the degree of $F$ implies that $\int_{\Sig_n}\Theta\wedge\o^{k+1}\wedge\ldots\wedge
\o^n=s$. 

We will need the following elementary lemma.
\begin{lem} \label{lem:exact}
Suppose $N^m$ is a closed oriented Riemannian manifold and let $\O$ be its volume form.
Given any open set $U$ of $N$ which is not dense in $N$, the form $\O$ is exact on $U$. Moreover,
given an open set $V$ compactly contained in $U$, we can find a closed $m$-form
$\O_1$ which agrees with $\O$ on $M\setminus U$ and such that $\O_1=0$ in $V$.
\end{lem}
\begin{pf}
Let $f$ be a smooth function which is equal to $1$ in $U$ and such that $\int_Nf\ d\O=0$.
Let $u$ be a solution of $\Delta u=f$ and let $\theta$ be the $(m-1)$-form $\theta=*du$. We then
have $d\theta=d*du=(\Delta u)\O$, so we have $d\theta=\O$ on $U$. 

To prove the last statement, we let $\zeta$ be a smooth cutoff function which is equal to $1$
in $V$ and has compact support in $U$. We then define $\O_1=\O-d(\zeta*du)$. We then
have $\O_1=0$ in $V$ and $\O_1$ differs from $\O$ by an exact form.
\end{pf}

We now restate the existence theorem.
\begin{thm}
For a manifold $M=\Sig_n$ as described above, there is a $\Lambda$-bounded,
partially regular, minimal $k$-slicing Moreover, if $k\leq j\leq n-1$ and $\Sig_j$ is regular, then 
$\int_{\Sig_j}\Theta\wedge\omega^{k+1}\wedge\ldots\wedge\omega^j=s$.
\end{thm}
\begin{pf} We begin with the $1$-form $\o^n$ and we integrate to get a map $u_n:\Sigma_n\to S^1$
so that $\o^n=du_n$. Let $t$ be a regular value of $u_n$ and consider the hypersurface $S_n=u_n^{-1}(t)$. Because the map $F$ has degree $s$ and we have normalized our forms in
$X\times T^{n-k}$ to have integral $1$, we see that $\int_{S_n}\Theta\wedge\o^{k+1}\wedge\ldots
\wedge\o^{n-1}=s$. Let $\Sig_{n-1}$ be a least volume cycle in $\Sig_n$ with the property that
$\int_{\Sig_n}\Theta\wedge\o^{k+1}\wedge\ldots\wedge\o^{n-1}=s$. The existence follows from
standard results of geometric measure theory.

Now suppose for $j\geq k$ we have constructed a partially regular minimal $j+1$ slicing 
with the property that there is a form $\Theta_{j+1}$ of compact support which is
cohomologous to $\Theta\wedge\o^{k+1}\wedge\ldots\wedge\o^{j+1}$ such that
$\int_{\Sig_{j+1}}\Theta_{j+1}=s$. Since the slicing is partially regular, we have that the
Hausdorff dimension of ${\cal S}_{j+1}$ is at most $j-2$, so it follows that the image
$F_j({\cal S}_{j+1})$ under the projection map $F_j:\Sig_n\to X\times T^{j-k}$ is a compact
set of Hausdorff dimension at most $j-2$. It follows from Lemma \ref{lem:exact} that the form
$\O\wedge dt^{k+1}\wedge\ldots\wedge dt^j$ is exact in a neighborhood $U$ of $F_j({\cal S}_{j+1})$,
given a neighborhood $V$ of $F_j({\cal S}_{j+1})$ which is compact in $U$ we can find a form 
$\O_j$ which is cohomologous to $\O\wedge dt^{k+1}\wedge\ldots\wedge dt^j$ and vanishes in $V$. Pulling back we see that $\Theta_j=F^*\O_j$ vanishes in a neighborhood of ${\cal S}_{j+1}$ and is cohomologous to $\Theta\wedge\o^{k+1}\wedge\ldots\wedge\o^j$. We let $u_{j+1}$ be the map gotten
by integrating $\o^{j+1}$ and consider its restriction to $\Sig_{j+1}$. Since $u_{j+1}$ is in $L^2$
with respect to the weight $\rho_{j+2}$, we see that $\rho_{j+1}=u_{j+1}\rho_{j+2}$ is integrable on
$\Sig_{j+1}$. It then follows from the coarea formula that we can find a regular value $t$
of $u_{j+1}$ in ${\cal R}_{j+1}$ so that the hypersurface $S_j\subset \Sig_{j+1}$ given by 
$S_j=u_{j+1}^{-1}(t)$ has finite $\rho_{j+1}$-weighted volume and satisfies $\int_{S_j}\Theta_j=s$.
We can then solve the minimization problem for the $\rho_{j+1}$-weighted volume among
integer multiplicity rectifiable currents $T$ with support in $\Sig_{j+1}$, with no boundary in
${\cal R}_{j+1}$, and with $T(\Theta_j)=s$. A minimizer for this problem gives us $\Sig_j$
and completes the inductive step for the existence.
\end{pf}
\begin{rem} The existence proof above does not specify the homology class of the minimizers
even if the minimizers are smooth since we are minimizing among cycles for which the
integral of $\Theta_j$ is fixed. In general there may be homology classes for which the integral
of $\Theta_j$ vanishes. We have chosen the class to do the minimization in order to avoid
a precise discussion of the homology of the singular spaces in which we are working. In
the following we give a more precise existence theorem which specifies the homology
classes and allows them to be general integral homology classes, possibly torsion classes.
\end{rem}

We now formulate and prove a more general existence theorem for minimal $k$ slicings.
In the theorem we let $[\Sig_n]$ denote the fundamental homology class in $H_n(\Sig_n,\mathbb Z)$
and, for a cohomology class $\a\in H^p(\Sig_n,\mathbb Z)$, we let $\a\cap [\Sig_n]$ denote
its Poincar\'e dual in $H_{n-p}(M,\mathbb Z)$.
\begin{thm} \label{thm:exst2} Let $\Sig_n$ be a smooth oriented manifold of dimension $n$
and let $k$ be an integer with $1\leq k\leq n-1$. Let $\a^1,\ldots,\a^{n-k}$ be cohomology classes
in $H^1(\Sig_n, \mathbb Z)$, and suppose that $\a^{n-k}\cap\a^{n-k-1}\cap\ldots\cap\a^1\cap[\Sig_n]
\neq 0$ in $H_n(\Sig_n,\mathbb Z)$. There exists a partially regular minimal $k$ slicing
with $\Sig_j$ representing the homology class $\a^{n-j}\cap\ldots\cap\a^1\cap [\Sig_n]$.
\end{thm}

\begin{pf}
Assume that we are given a partially regular $\Lambda$-bounded minimal $(k+1)$-slicing 
which represents $\a_1,\ldots,\a_{n-k-1}$. We thus have the weight 
function $\rho_{k+1}$ defined on $\Sig_{k+1}$ which we use to produce $\Sig_k$. From the
partial regularity the singular set ${\cal S}_{k+1}$ of $\Sig_{k+1}$ has Hausdorff dimension
at most $k-2$. 

We consider the class of integer multiplicity rectifiable currents which are relative cycles in
$H_k(\Sig_n,{\cal S}_{k+1},\mathbb Z)$; that is, for any $k-1$ form $\theta$ of compact support in
$\Sig_{k+1}\setminus{\cal S}_{k+1}$ we have $T(d\theta)=0$. Because the set ${\cal S}_{k+1}$
has zero $k-1$ dimensional Hausdorff measure we have $H_k(\Sig_n,{\mathbb Z})=
H_k(\Sig_n,{\cal S}_{k+1},\mathbb Z)$. This follows because a current which is a relative
cycle $T$ in $\Sig_n\setminus{\cal S}_{k+1}$ is also a cycle in $\Sig_n$ since $\p T$
is zero since it is unchanged by adding a set of $k-1$ measure zero.

We use $\rho_{k+1}$ weighted volume to set up a minimization problem. We consider the
class of relative cycles $T$ with support contained in $\Sig_{k+1}$ which have finite weighted mass;
that is, $T=(S_k,\Theta,\xi)$ where $S_k$ is 
a countably $k$-rectifiable set, $\Theta$ a $\mu_k$-measurable integer valued function on $S_k$,
and $\xi$ a $\mu_k$-measurable map from $S_k$ to $\wedge^k\R^N$ such that $\xi(x)$ is a
unit simple vector for $\mu_k$ a.e. $x\in S_k$. Such a $k$-current $T_k$ is $\rho_{k+1}$-finite if
\[ Vol_{\rho_{k+1}}(T_k)\equiv\int_{S_k}\rho_{k+1}|\Theta|\ d\mu_k<\infty.
\]

Since we have already constructed $\Sig_{k+1}$ so that it is $\Lambda$-bounded we have
\[ \int_{\Sig_{k+1}}\rho_{k+1}\ d\mu_{k+1}\leq \Lambda.
\]
Now we can find a smooth closed hypersurface $H_k$ which is Poincar\'e dual to $\a_k$, and 
we may perturb
it and use the coarea formula in a standard way to arrange that $\bar{\Sig}_k\equiv\Sig_{k+1}\cap H_k$
is a smooth embedded submanifold away from  ${\cal S}_{k+1}$ and
\[ \int_{\bar{\Sig}_k}\rho_{k+1}\ d\mu_k\leq c.
\]
In particular the associated current $\bar{T}_k\equiv(\bar{\Sig}_k,1,\bar{\xi})$ (where $\bar{\xi}$ is the oriented unit tangent plane of $\bar{\Sig}_k$) is $\rho_{k+1}$-finite and is a competitor in our
variational problem. 

The standard theory of integral currents now allows us to construct a minimizer for our variational
problem which gives us the next slice $\Sig_k$ which could be disconnected and with integer
multiplicity. Thus $\Sig_k$ represents the homology class 
$\a^{n-k}\cap\ldots\cap\a^1\cap [\Sig_n]$. This completes the proof of Theorem \ref{thm:exst2}.
\end{pf}

\bigskip
\section{\bf Application to scalar curvature problems}
In this section we prove two theorems for manifolds with positive scalar curvature. The first of these
is for compact manifolds and the second is the Positive Mass Theorem for asymptotically flat 
manifolds. Our first theorem which we will need to prove the Positive Mass Theorem is the following.
\begin{thm} \label{thm:psc0} Let $M_1$ be any closed oriented $n$-manifold. The manifold
$M=M_1\#T^n$ does not have a metric of positive scalar curvature.
\end{thm}
\begin{pf} Such a manifold $M$ has admits a map $F:M\to T^n$ of degree $1$, and so by
Theorem \ref{thm:exst} there exists a closed minimal $1$-slicing of $M$ in contradiction to
Theorem \ref{thm:12slicing}.

\end{pf}
We also prove the following more general theorem.
\begin{thm}
\label{thm:psc1} Assume that $M$ is a compact oriented $n$-manifold with a metric of positive
scalar curvature. If $\a_1,\ldots,\a_{n-2}$ are classes in $H^1(M,\Z)$ with the property that the
class $\s_2$ given by
$\s_2=\a_{n-2}\cap\a_{n-3}\cap\ldots\a_1\cap[M]\in H_2(M,\Z)$ is nonzero, then the class 
$\s_2$ can be represented by a sum of smooth two spheres.
If $\a_{n-1}$ is any class in $H^1(M,\Z)$, then we must have $\a_{n-1}\cap\s_2=0$. In particular,
if $M$ has classes $\a_1,\ldots,\a_{n-1}$ with $\a_{n-1}\cap\ldots\cap\a_1\cap[M]\neq 0$,
then $M$ cannot carry a metric of positive scalar curvature.
\end{thm}
\begin{pf} By the existence and regularity results of Sections 3 and 4, there is a minimal
$2$-slicing so that $\Sig_2\in\s_2$ is regular and satisfies the eigenvalue bound of Theorem \ref{thm:eval}. Choosing $\varphi=1$ on any given component of $\Sig_2$ and applying the Gauss-Bonnet
theorem we see that each component must be topologically $S^2$. 

In particular it follows that for any other $\a_{n-1}\in H^1(M,\Z)$ we have that $\a_{n-1}\cap\s_2$
is a class in $H_1(\Sig_2,\Z)$, and therefore is zero. 
\end{pf}

We now prove a Riemannian version of the positive mass theorem. Assume that $M$ is a 
complete manifold with the property that there is a compact subset $K\subset M$ such that
$M\sim K$ is a union of a finite number of connected components each of which is an
asymptotically flat end. This means that each of the components is diffeomorphic to the
exterior of a compact set in $\R^n$ and admits asympototically flat coordinates $x^1,\ldots,x^n$
in which the metric $g_{ij}$ satisfies
\begin{equation}
\label{eqn:af}
g_{ij}=\d_{ij}+O(|x|^{-p}),\ |x||\partial g_{ij}|+|x|^2|\partial^2g_{ij}|=O(|x|^{-p}),\ |R|=O(|x|^{-q})
\end{equation}
where $p>(n-2)/2$ and $q>n$. Under these assumptions
the ADM mass is well defined by the formula (see \cite{sc} for the $n$ dimensional case)
\[ m=\frac{1}{4(n-1)\o_{n-1}}\lim_{\s\to\infty}\int_{S_\s}\sum_{i,j}(g_{ij.i}-g_{ii,j})\nu_j\ d\xi(\s)
\]
where $S_\s$ is the euclidean sphere in the $x$ coordinates, $\o_{n-1}=Vol(S^{n-1}(1))$, and 
the unit normal and volume integral are with respect to the euclidean metric. We may now state
the Positive Mass Theorem.
\begin{thm}
\label{thm:psc2} Assume that $M$ is an asymptotically flat manifold with $R\geq 0$. For each end
it is true that the ADM mass is nonnegative. Furthermore, if any of the masses is zero, then $M$ is isometric to $\R^n$.
\end{thm}
\begin{pf} The theorem can be reduced to the case when there is a single end by capping off
the other ends keeping the scalar curvature nonnegative. We will show only that $m\geq 0$,
and the equality statement can be derived from this (see \cite{sy2}). We will reduce the proof
to the compact case using results of \cite{sy3} and an observation of J. Lohkamp. 
\begin{prop}
If the mass of $M$ is negative, there is a metric of nonnegative scalar curvature on $M$ which
is euclidean outside a compact set. This produces a metric of positive scalar curvature on a manifold
$\hat{M}$ which is gotten by replacing a ball in $T^n$ by the interior of a large ball in $M$. 
\end{prop}
\begin{pf} Results of \cite{sy3} and \cite{sc} imply that if $m<0$ we can construct a new metric
on $M$ with nonnegative scalar curvature, negative mass, and which is conformally flat and
scalar flat near infinity. In particular, we have $g=u^{4/(n-2}\d$ near infinity where $u$ is a
euclidean harmonic function which is asymptotic to $1$. Thus $u$ has the expansion
\[ u(x)=1+\frac{m}{|x|^{n-2}}+O(|x|^{1-n})
\]
where $m$ is the mass. Now we use an observation of Lohkamp \cite{lohkamp}.  Since $m<0$, we can choose
$0<\e_2<\e_1$ and $\s$ sufficiently large so that we have $u(x)<1-\e_1$ for $|x|=\s$ and
$u(x)>1-\e_2$ for $|x|\geq 2\s$. If we define $v(x)=u(x)$ for $|x|\leq \s$ and $v(x)=\min\{1-\e_2,u(x)\}$
for $|x|>\s$, then we see that $v(x)$ is weakly superharmonic for $|x|\geq \s$, so may be 
approximated by
a smooth superharmonic function with $v(x)=u(x)$ for $|x|\leq \s$ and $v(x)=1-\e_2$ for $|x|$ 
sufficiently large. The metric
which agrees with the original inside $S_\s$ and is given by $v^{4/(n-2)}\d$ outside then has
nonnegative scalar curvature and is euclidean near infinity.

By extending this metric periodically we then produce a metric on $\hat{M}$ with nonnegative scalar
curvature which is not Ricci flat. Therefore the metric can be perturbed to have positive 
scalar curvature.
\end{pf}
Using this result the theorem follows from Theorem \ref{thm:psc1} since the standard $1$-forms
on $T^n$ can be pulled back to $\hat{M}$ to produce the $\a_1,\ldots,\a_{n-1}$ of that
theorem. This completes the proof of Theorem \ref{thm:psc2}.
\end{pf}

\end{document}